\documentclass[12pt]{amsart}
\usepackage{graphicx} 

\usepackage[T1]{fontenc}
\usepackage[utf8]{inputenc}
\usepackage[english]{babel}
\usepackage{amsmath,amsthm,amssymb,latexsym,amsfonts, mathtools}
\usepackage{mathrsfs} 
\usepackage{tikz-cd}
\usepackage{framed}
\usepackage{listings}
\usepackage{shuffle}
\usepackage{faktor} 
\usepackage{nicefrac}
\usepackage{color,soul}
\usepackage{setspace}
\usepackage{graphicx}
\usepackage{faktor}
\usepackage[a4paper, margin=1in]{geometry}
\usepackage[alphabetic]{amsrefs}
\usepackage[colorlinks,linkcolor=black,anchorcolor=blue,citecolor=blue]{hyperref}

\usepackage{xpatch}
\makeatletter   
\xpatchcmd{\@tocline}
{\hfil\hbox to\@pnumwidth{\@tocpagenum{#7}}\par}
{\ifnum#1<0\hfill\else\dotfill\fi\hbox to\@pnumwidth{\@tocpagenum{#7}}\par}
{}{}
\makeatother

\usepackage{etoolbox}
\makeatletter
\patchcmd\@thm
  {\let\thm@indent\indent}{\let\thm@indent\noindent}%
  {}{}
\makeatother

\makeatletter
\patchcmd{\@maketitle}
  {\ifx\@empty\@dedicatory}
  {\ifx\@empty\@date \else {\vskip3ex \centering\footnotesize\@date\par\vskip1ex}\fi
   \ifx\@empty\@dedicatory}
  {}{}
\patchcmd{\@adminfootnotes}
  {\ifx\@empty\@date\else \@footnotetext{\@setdate}\fi}
  {}{}{}
\makeatother

\makeatletter
\patchcmd{\@thm}{\thm@headfont{\scshape}}{\thm@headfont{\scshape\normalfont\bfseries}}{}{}
\patchcmd{\@thm}{\thm@notefont{\fontseries\mddefault\upshape}}{}{}{}
\makeatother
\usepackage{subfiles}

\DeclareFontFamily{U}{wncy}{}
\DeclareFontShape{U}{wncy}{m}{n}{<->wncyr10}{}
\DeclareSymbolFont{mcy}{U}{wncy}{m}{n}
\DeclareMathSymbol{\sh}{\mathord}{mcy}{"78}

\renewcommand{\epsilon}{\varepsilon}
\renewcommand{\theta}{\vartheta}
\renewcommand{\rho}{\varrho}
\let\temp\phi
\let\phi\varphi
\let\varphi\temp
\newcommand{\insubfile}[1]{\ifx\@onlypreamble\@notprerr\else#1\fi}
\newcommand{\f}[1]{\mathfrak{#1}}

\newcommand{\Z}{\mathbb{Z}}
\newcommand{\Q}{\mathbb{Q}}
\newcommand{\R}{\mathbb{R}}
\newcommand{\C}{\mathbb{C}}
\newcommand{\F}{\mathcal{F}}
\newcommand{\A}{\mathcal{A}}

\renewcommand{\S}{\mathcal{S}}
\newcommand{\modulo}[1]{\;(\textnormal{mod}\; #1)}

\theoremstyle{definition}
\newtheorem{definition}{Definition}
\theoremstyle{plain}
\newtheorem{proposition}[definition]{Proposition}
\theoremstyle{plain}
\newtheorem{theorem}[definition]{Theorem}
\theoremstyle{plain}
\newtheorem{lemma}[definition]{Lemma}
\theoremstyle{plain}
\newtheorem{corollary}[definition]{Corollary}
\theoremstyle{plain}

\theoremstyle{definition}
\newtheorem{lemma*}{Lemma}
\theoremstyle{definition}
\newtheorem{remark}[definition]{Remark}
\theoremstyle{definition}

\AtBeginDocument{%
   \def\MR#1{}
}

\title{An explicit study of a family of cellular integrals}
\author{Riccardo Tosi}

\address[Riccardo Tosi]{ University of Duisburg-Essen, Fakultät für Mathematik, Thea-Leymann-Str. 9, 45127 Essen, Germany.}

\email[]{riccardo.tosi@uni-due.de}

\begin{document}

    \begin{abstract}
        
    We express a family of basic cellular integrals over moduli spaces of curves explicitly in terms of multiple zeta values, answering a question of Brown. Moreover, we study a priori the weights appearing in these integrals and find a relation that expresses the odd-dimensional integrals in terms of the even-dimensional ones. We also sketch an explanation of this relation in the spirit of Grothendieck's Period Conjecture. 

    \end{abstract}
    \noindent
     
    \maketitle

    
    \tableofcontents

    \section*{Introduction}
     
    The values of the Riemann zeta function at even positive integers are rational multiples of powers of $\pi$, hence transcendental. Little is known about the arithmetic nature of $\zeta(s)$ for $s\ge 3$ odd: among the few available results, there is the irrationality of $\zeta(3)$ \cite{Apery-Irrationalité_de_zeta2_et_zeta3} and the infinite dimension of the $\Q$-vector space generated by odd zeta values \cite{Ball-Rivoa-Irrationalité_d'une_infinité_de_valuers_de_la_fonction_zeta_aux_entiers_impairs}. Typical methods to address irrationality questions involve constructing sequences of $\Q$-linear combinations of odd zeta values which tend to zero relatively fast compared to the denominators of the coefficients.\\
    A geometric input to this problem comes from the periods of moduli spaces of curves of genus zero, which are known to be $\Q(2\pi i)$-linear combinations of multiple zeta values \cites{Brown-Multiple_zeta_values_and_periods_of_moduli_spaces_of_curves, Brown-Mixed_Tate_motives_over_Z}. Thus, integrals of algebraic differential forms over these varieties may serve as a natural source of linear forms necessary for irrationality proofs of (multiple) zeta values. Brown \cite{Brown-Irrationality_proofs_for_zeta_values_moduli_spaces_and_dinner_parties} has isolated certain promising families of these integrals, called \emph{cellular integrals}, which present large groups of symmetries and good analytic properties. A study of one of these families \cite{Brown-Zudilin-Cellular_rational_approximations_to_zeta5} has led, for example, to the best rational approximations of $\zeta(5)$ known to date. \\
    Let $\f{M}_{0,l+3}$ denote the moduli space of smooth curves of genus zero with $(l+3)$-marked points. Brown's cellular integrals take the shape, for $m_1,\dots, m_r\in \Z$,
    \[
        \int_{\Delta} f_1^{m_1}\dots f_r^{m_r} \omega,
    \]
    where $\Delta$ is a relative homology class, $f_1,\dots, f_r$ are some elementary algebraic functions on $\f{M}_{0,l+3}$ and $\omega$ is a volume form on an appropriate partial compactification of $\f{M}_{0,n}$. Panzer \cite{Panzer-Algortihms_for_the_symbolic_integration_of_hyperlogarithms_with_applications_to_Feynman_integrals} has developed an algorithm to compute these integrals that has been implemented in the program \texttt{HyperInt}. However, the task of computing infinite families of cellular integrals, as required by irrationality proofs, remains rather difficult. \\  
    In this article, we provide an explicit description in terms of multiple zeta values of the simplest possible cellular integrals, corresponding in the above notation to the choice $f_1=\dots=f_r=1$. In suitable coordinates, for $l\ge 2$ these integrals are given by
    \[
        \xi_l=\int_{[0,1]^l}\frac{1}{(1-x_1x_2)(1-x_2x_3) \dots (1-x_{l-1}x_l)} \, dx_1\dots dx_l.
    \]
    Our main result, which answers a question of Brown \cite[Remark 8.7]{Brown-Multiple_zeta_values_and_periods_of_moduli_spaces_of_curves}, reads as follows.
    \begin{theorem}
    \label{theo: main}
        For all $m\ge 1$, the integrals $\xi_{2m}$ fit into the following generating series:
        \[
            1+\sum_{n=1}^\infty \xi_{2n}t^n= \frac{\arcsin \left( \pi \sqrt{t}\right)}{\pi \sqrt{t}}.
        \]
        Moreover, 
        \[
            \xi_{2m+1}=\sum_{h=0}^{m} \xi_{2h}\xi_{2m-2h}.
        \]
    \end{theorem}
    \noindent
    In particular, each $\xi_{l}$ is a rational multiple of $\pi^l$ for $l$ even and of $\pi^{l-1}$ for $l$ odd. For $l=2m$, we deduce the following formula:
    \[
        \xi_{2m}= \frac{(2m)!}{4^m (m!)^2(2m+1)}\, \pi^{2m}.
    \]
    It is worth pointing out that for $l$ even the integral $\xi_l$ is a linear combination of multiple zeta values of weight exactly $l$, while for $l$ odd the weight is exactly $l-1$. \\
    \noindent
    The paper is organized as follows. In the first section, we study the behaviour of the weights appearing in $\xi_l$. For $l$ even, the vanishing of weights lower than $l$ is a consequence of some elementary facts about mixed Hodge Tate structures. For $l$ odd, the absence of multiple zeta values of weight $l$ is due to the integrand of $\xi_l$ being exact as an algebraic form. Apart from exhibiting a primitive thereof, we also sketch how to justify the relation between $\xi_{2m+1}$ and $\xi_{2m'}$ with $m'<m$ in terms of Grothendieck's period conjecture. \\
    In the second section, we prove Theorem~\ref{theo: main} by means of an inductive application of Panzer's algorithm \cite{Panzer-Algortihms_for_the_symbolic_integration_of_hyperlogarithms_with_applications_to_Feynman_integrals}. We compute explicitly a few examples in low dimension, then study the combinatorics of a specific recurrence sequence to address the general case. The arguments do not strictly need the results of the first section, but studying the weights of $\xi_l$ a priori is of independent interest for irrationality proofs. \\
    The program proposed in \cite{Brown-Irrationality_proofs_for_zeta_values_moduli_spaces_and_dinner_parties} to tackle the irrationality of zeta values through cellular integrals still presents several challenges, as it is difficult to predict which weights appear in a given integral. We hope that understanding the most basic family of these integrals very explicitly may help to make the general case more treatable. It would be interesting to see which cellular integrals can be studied by reducing their computation to the one of the sequence $\xi_l$ of the present paper. \\

    \noindent
    \textbf{Acknowledgements.} The author expresses his heartfelt gratitude to Johannes Sprang for the careful supervision of this project and his constant encouragement. The author is also indebted to Clément Dupont, without whose suggestions the first section of this paper would not have been written, and to Erik Panzer, for pointing out a mistake in a previous version of this article and proposing further improvements. Yasuo Ohno and Wadim Zudilin's comments have also improved this work. The support of the DFG Research Training Group 2553 is gratefully acknowledged.

    \section{Study of the weights}
    \subsection{Moduli spaces of curves of genus zero}
        
    We recall two different descriptions of the moduli spaces of smooth and stable curves of genus zero with marked points. Let $l\ge 1$ be an integer and define $n=l+3$. Let $V$ be an $(l+1)$-dimensional $\Q$-vector space and denote by $V^*$ its dual. Fix a basis $y_1,\dots, y_{l+1}$ of $V^*$. Consider the set of lines in $V^*$ given by
    \[
        \A_l=\{\langle y_i \rangle \mid i=1,\dots, l \} \cup \{ \langle y_i-y_j\rangle \mid i,j=1,\dots, l+1,\; i\neq j\}.
    \]
    Each of these lines $H\in \A$ defines a hyperplane $H^\perp$ in $V$ by considering the kernel of any non-zero element of $H$. The moduli space $\f{M}_{0,n}$ of smooth projective curves over $\Q$ of genus $0$ with $l+3$ marked points is isomorphic to the projective complement of the hyperplanes associated with $\A_l$, that is,
    \[
        \f{M}_{0,n}\cong \mathbb{P}(V)\setminus \bigcup_{H\in \A} \mathbb{P}(H^\perp).
    \]
    Explicitly, in affine coordinates $t_i=\frac{y_i}{y_{l+1}}$, usually called \emph{simplicial coordinates}, we have
    \[
        \f{M}_{0,n} \cong \text{Spec}\; \Q[t_1,\dots,t_l]\left[\frac{1}{t_i},\frac{1}{1-t_i},\frac{1}{t_i-t_j} \;\middle\vert\; i,j=1,\dots l, \; i\neq j \right].
    \]
    De Concini and Procesi \cite{De_Concini-Procesi-Wonderful_models_of_subspace_arrangements} described a compactification of $\f{M}_{0,n}$ by a simple normal crossing divisor explicitly in terms of the combinatorics of the arrangement $\A_{l}$. This compactification is isomorphic to the moduli space $\overline{\f{M}}_{0,n}$ of stable curves of genus $0$ with $l+3$ marked points. We summarize a few results from \cite{De_Concini-Procesi-Wonderful_models_of_subspace_arrangements}.\\
    Let $L(\A_l)$ be the intersection lattice of $\A_l$, that is, the set of all subspaces $X$ of $V^*$ which are sums of elements of $\A_l$, which we regard as a partially ordered set with respect to standard inclusion. Given a subset $\lambda\subseteq \{ 0,1,\dots, l+1\}$ with $\#\lambda \ge 2$, we define the following element of $L(\A_l)$:
    \[
        Y_\lambda=
        \begin{cases}
            \langle y_i \mid i\in \lambda\setminus \{0\}\rangle & \text{if $0\in \lambda$}; \\
            \langle y_i-y_j\mid i,j\in \lambda, \; i\neq j \rangle & \text{if $0\not\in \lambda$}.
        \end{cases}
    \]
    We call these elements of $L(\A_l)$ \emph{irreducible} and denote their set by $\F_l$. Every $X\in L(\A_l)$, $X\neq 0$, can be written uniquely as the direct sum of finitely many irreducible elements of $L(\A_l)$. Thus, $L(A_l)\setminus \{0\}$ is isomorphic, as a partially ordered set, to the set of partitions $\{\lambda_1,\dots, \lambda_k\}$ of subsets $\Lambda\subseteq \{0,\dots, l+1\}$ with $\#\lambda_i \ge 2$, where the partial order is given by $\{ \lambda_1,\dots ,\lambda_k\} \le \{\lambda'_1,\dots, \lambda'_{k'}\}$ if and only if there are $1\le i_1<\dots < i_h \le k'$ such that $\bigcup_{j=1}^k \lambda_i = \bigcup_{j=1}^h \lambda'_{i_j}$ and $\{ \lambda_1,\dots, \lambda_k\}$ is a refinement of the partition $\{\lambda'_{i_1},\dots, \lambda'_{i_h}\}$. \\
    De Concini and Procesi construct a compactification of $\f{M}_{0,n}$, which we denote by $\overline{\f{M}}_{0,n}$, by iteratively blowing up $\mathbb{P}(V)$ along the successive strict transforms of the $Y_\lambda$'s in a precise order. They develop a general method to compactify complements of arrangements of linear subspaces by blowing up along appropriately defined irreducible elements, but, for our purposes, we will describe their results only for $\f{M}_{0,n}$.\\
    For all $X\in \F_l$ there is a well defined morphism $V\setminus X^\perp \to \mathbb{P}(V/X^\perp)$. This yields a morphism $\f{M}_{0,n}\to \mathbb{P}(V/X^\perp)$. We then take $\overline{\f{M}}_{0,n}$ to be the closure of the graph of the product of these morphisms inside $\mathbb{P}(V) \times \prod_{X\in \F_l} \mathbb{P}(V/X^\perp)$. This is an irreducible projective variety of which $\f{M}_{0,n}$ makes up an open subscheme. In order to work with explicit local charts and describe the boundary divisor, we need to introduce some combinatorial notions attached to the arrangement $\A_l$.\\
    A subset $\S\subseteq L(\A_l)\setminus\{0\}$ is called \emph{$\F_l$-nested} if every $X\in \S$ is irreducible and for all $X_1,\dots, X_k\in \S$ pairwise non-comparable we have $\sum_{i=1}^kX_i$ is not irreducible. The latter condition is equivalent to requiring that for all $X=Y_\lambda, X'=Y_{\lambda'}\in \S$ with $\lambda\not\subseteq \lambda'$, $\lambda'\not\subseteq \lambda$ we have $\lambda\cap \lambda'=\emptyset$. \\
    Given a $\F_l$-nested set $\S$, an \emph{adapted basis} for $\S$ is a function $\beta\colon \S \to V^* $ such that $X\in \S$ the set $\{ \beta(Y) \mid Y\in \S, Y\subseteq X\}$ is a basis for $X$. Adapted bases always exist and can be chosen so that $\beta(X)$ belongs to an element of $\A_l$ for all $X\in \S$. \\
    Let $\S$ be a $\F_l$-nested set which is maximal with respect to inclusion. It can be seen that $\#\S=l+1$ and that every $\F_l$-nested set can be completed to a maximal one. Moreover, for all $X\in \S$, $X\neq V^*$, the set of $Y\in \S$ that strictly contain $X$ is linearly ordered, hence it has a minimum $X^+$. \\
    Fix an adapted basis $\beta$ for $\S$. Consider the morphism
    \[
        \rho \colon  \mathbb{A}^{l+1}=\text{Spec}\, \Q[u_X\mid X\in \S]\longrightarrow \mathbb{A}^{l+1}=\text{Spec}\, \Q[\beta(X)\mid X\in \S]
    \]
    which, at the level of regular functions, is defined by
    \[
        \beta(X)\longmapsto \prod_{Y\in \S, X\subseteq Y} u_Y.
    \]
    The map $\rho$ restricts to an isomorphism between the open subsets defined by removing the hyperplanes $u_X=0$ in the source and $\beta(X)=0$ in the target for all $X\in \S$. Its inverse is given by
    \[
        u_X\longmapsto \frac{\beta(X)}{\beta(X^+)},
    \]
    formally setting $\beta((V^*)^+)=1$. \\
    Let $H\in \A$ and $x\in H$, $x\neq 0$. There is a minimal element $p_\S(x)\in \S$ to which $x$ belongs. Then $x=\beta(p_\S(x)) P^\S_x$, where the image of $P^\S_x$ under $\rho$ is a polynomial that depends only on the coordinates $u_Y$ for $Y\in \S$, $Y\subsetneq p_\S(x)$ and does not vanish at $0$. \\
    We define the affine variety
    \[ 
        U_{\S}=\text{Spec}\, \Q[u_X\mid X\in \S, X\neq V^*]\left[ \left(P^\S_H\right)^{-1} \,\middle\vert H\in \A \right],
    \]
    where $P^\S_H=P^\S_x$ for a choice of $x\in H$, $x\neq 0$. This choice does not affect the definition of $U_\S$. The map $\rho$ induces an isomorphism
    \[
        U_\S \setminus \bigcup_{X\in \S} \{ u_X=0\} \cong \f{M}_{0,n}.
    \]
    Each $U_\S$ embeds naturally as an open subset of $\overline{\f{M}}_{0,n}$.
    Moreover, the local charts $U_\S$ for $\S$ ranging among all maximal $\F_l$-nested sets of $L(\A_l)$ cover the projective variety $\overline{\f{M}}_{0,n}$.
    \begin{theorem}[{\cite[Proposition 1.5]{De_Concini-Procesi-Wonderful_models_of_subspace_arrangements}}]
         $ $
        \begin{enumerate}
            \item The variety $\overline{\f{M}}_{0,n}$ is smooth and the complement of $\f{M}_{0,n}$ therein is a simple normal crossings divisor. 
            \item The irreducible components of $\overline{\f{M}}_{0,n}\setminus \f{M}_{0,n}$ are in bijection with the elements of $\F_l \setminus \{V^*\}$. The irreducible divisor $D_X$ corresponding to $X\in \F_l$ is the closure of $\{u_X=0\}\subseteq U_\S$ for any maximal $\F_l$-nested set $\S$ containing $X$. 
            \item The irreducible divisors $D_{X_1},\dots, D_{X_k}$  corresponding to $X_1,\dots, X_k\in \F_l\setminus \{V^*\}$ have non-empty intersection if and only if $\{X_1,\dots, X_k\}$ is $\F_l$-nested. Moreover, this intersection is smooth. 
        \end{enumerate}
    \end{theorem}
    \noindent 
    The choice $\S= \bigcup_{i=1}^{l+1} \{ \langle y_1,\dots, y_i \rangle\}$ with adapted basis $\beta(\langle y_1,\dots, y_i \rangle)=y_i$ yields the local chart
    \[
        U_\S=\text{Spec}\, \Q[x_1,\dots, x_l]\left[ \frac{1}{1-x_i\dots x_j} \;\middle\vert\; 1\le i \le j \le l \right],
    \]
    where we have written $x_i$ for the coordinate corresponding to $\langle y_1,\dots , y_i\rangle$. These are usually referred to as \emph{cubical coordinates}. Note that, in the real points of $\f{M}_{0,n}$, the simplex $0<t_1<\dots < t_l <1$ in simplicial coordinates corresponds to the cube $\prod_{i=1}^l \{ 0<x_i<1\}$ in cubical coordinates. \\
    The boundary divisors of a De Concini-Procesi compactification are themselves isomorphic to products of compactifications of complements of hyperplane arrangements of smaller dimension. For our purposes, we only recall that for all $X\in \F(\A_l)$ the associated divisor $D_X$ of $\overline{\f{M}}_{0,n}$ decomposes as
    \[
        D_X\cong  \overline{\f{M}}_{0,n_1} \times  \overline{\f{M}}_{0,n_2}
    \]
    for some $n_1,n_2\in\{0,\dots, n-1\}$ such that $n_1+n_2=n-1$.\\
    For the purpose of computing integrals, we describe some canonical relative homology classes of $\f{M}_{0,n}$. Let $\delta$ be a permutation of the set $\{0,1,\dots, l+1\}$ and write $<_\delta$ for the linear order induced by $\delta$ on $\{0,\dots, l+1\}$, that is, $i<_\delta j$ if and only if $\delta^{-1}(i)<\delta^{-1}(j)$. Suppose that $0<_\delta l+1$. Setting $t_0=0$ and $t_{l+1}=1$, there is a corresponding connected component of $\f{M}_{0,n}(\R)$, namely
    \[
        X_{n,\delta} = \{ (t_1,\dots, t_l)\in \R^l \mid 
        t_0<_\delta t_1<_\delta \dots <_\delta t_{l+1} \}.
    \]
    Moreover, all connected components of $\f{M}_{0,n}$ arise in this way. Thus, these components are in bijection with the linear orders of $\{0,\dots, l+1\}$ modulo the equivalence relation that identifies opposite orders. \\
    Let $\overline{X}_{n,\delta}$ be the closure of $X_{n,\delta}$ in $\overline{\f{M}}_{0,n}(\R)$. For $X\in \F_l$, we say that the irreducible boundary divisor $D_X$ of $\overline{\f{M}}_{0,n}$ is \emph{at finite distance from $\overline{X}_{n,\delta}$} if $\overline{X}_{n,\delta}\cap D_X(\R)\neq \emptyset$. It is easily checked that these divisors are precisely the ones corresponding to $Y_\lambda \in \F_l$ where $\lambda\subsetneq \{0,\dots, l+1\}$ is a segment with respect to $<_\delta$. Thus, $\overline{X}_{n,\delta}$ defines a singular homology class of degree $l$ of $\overline{\f{M}}_{0,n}$ relative to the irreducible boundary divisors at finite distance from it. \\

    The description of $\f{M}_{0,n}$ given so far only depends on the combinatoric of the underlying arrangement of hyperplanes and can therefore be generalized to other arrangements. However, it is possible to find more symmetric local coordinates on $\f{M}_{0,n}$, which exhibit a richer structure of its automorphism group. We follow the exposition of \cite{Brown-Multiple_zeta_values_and_periods_of_moduli_spaces_of_curves}. \\
    Let $(\mathbb{P}^1)^n_*$ be the product of $n$ copies of $\mathbb{P}^1$ without the big diagonal. This means that, giving coordinates $z_1,\dots, z_n$ to each copy of $\mathbb{P}^1$, one removes from $(\mathbb{P}^1)^n$ the closed subschemes given by $z_i=z_j$ for $i\neq j$. If we let $\textnormal{PSL}_2$ act diagonally on $(\mathbb{P}^1)^n_*$, we have 
    \[
        \f{M}_{0,n}\cong (\mathbb{P}^1)^n_* / \textnormal{PSL}_2.
    \]
    The regular functions of $\f{M}_{0,n}$ are the $\textnormal{PSL}_2$-invariant regular functions on $(\mathbb{P}^1)^n_*$. The latter are generated by \emph{cross-ratios}, that is, by functions of the form
    \[
        [ij|hk]=\frac{(z_i-z_h)(z_j-z_k)}{(z_i-z_k)(z_j-z_h)}
    \]
    for distinct $i,j,h,k\in \{1,\dots, n\}$. Since the action of $\textnormal{PSL}_2$ on $\mathbb{P}^1$ is triply transitive, we may recover the above simplicial coordinates by putting $z_1=1$, $z_2=\infty$ and $z_3=0$, while setting $t_1=z_4,\dots, t_l=z_n$.\\
    The action of the symmetric group $\textnormal{Sym}_n$ on $n$ letters over $(\mathbb{P}^1)^n_*$ by permuting coordinates transfers to $\f{M}_{0,n}$. This group acts transitively on the connected components of $\f{M}_{0,n}(\R)$ with stabilizer given by the dihedral groups $D_{2n}$ of symmetries of an $n$-gon.  \\
    These automorphisms extend to the compactification $\overline{\f{M}}_{0,n}$ and induce an action of $\textnormal{Sym}_n$ on the set of irreducible boundary divisors of $\overline{\f{M}}_{0,n}$. To describe this action, let us identify $\{0,\dots, l+1\}$ with the set $\{z_1,\dots, z_n\} \setminus \{z_2\}$ via $i \mapsto z_{i+3}$ if $i\neq l+1$ and $l+1\mapsto z_1$, in accordance with the above convention. An irreducible boundary divisor of $\overline{\f{M}}_{0,n}$ corresponds to a proper subset $\lambda \subsetneq \{0,\dots, l+1\}$ with $\#\lambda \ge 2$. This uniquely determines a partition of $\{z_1,\dots, z_n\}$ into two disjoint subsets $S_1,S_2$ with $\# S_1,\#S_2\ge 2$: to recover the original subset of $\{0,\dots,l+1\}$ it suffices to consider the set between $S_1$ and $S_2$ that does not contain $z_2$. The group $\text{Sym}_n$ acts naturally on these partitions, hence on the boundary divisors of $\overline{\f{M}}_{0,n} $, by permuting the $z_i$'s. \\
    Fix $\delta\in \text{Sym}(\{0,\dots, l+1\})$. We may identify the set $\{z_1,\dots, z_n\}$ with the edges of an $n$-gon by choosing an edge for $z_2$, followed by $z_{\delta^{-1}(3)+3}$ and so on. The boundary divisors at finite distance from $\overline{X}_{n,\delta}$ are given by proper segments of $\{0,\dots, l+1\}$ with respect to $<_\delta$, which correspond therefore to the diagonals of the $n$-gon. The dihedral group $D_{2n}\subseteq \textnormal{Sym}_n$ acts on these diagonals via the symmetries of the $n$-gon; in particular, it permutes the divisors at finite distance from $\overline{X}_{n,\delta}$.
    \begin{remark}
        The moduli description of $\f{M}_{0,n}$ makes the action of the whole group $\text{Sym}_n$ visible, while the interpretation via hyperplane arrangements gives access to fewer automorphisms, namely $\text{Sym}_{n-1}$. We will distinguish these two descriptions and rely only on the latter when possible, as this allows for generalizations to other hyperplane arrangements.
    \end{remark}
    \noindent
    Consider the following algebraic differential $n$-form on $(\mathbb{P}^1)^n_*$:
    \[
        \widetilde{\omega}_n= \prod_{i=1}^n \frac{1}{z_i-z_{i+2}} \, dz_1 \dots dz_n,
    \]
    where the indices in the product are considered modulo $n$. This form is invariant under the actions of both $\textnormal{PGL}_2$ and the dihedral group $D_{2n}$. \\
    Let $\nu$ be a non-zero algebraic invariant $3$-form on $\textnormal{PSL}_2$, which is unique up to a rational multiple. Since the quotient map $(\mathbb{P}^1)^n_*\to \f{M}_{0,n}$ is a trivial $\textnormal{PSL}_2$-torsor, we have $(\mathbb{P}^1)^n_* \cong \f{M}_{0,n} \times \textnormal{PSL}_2$. The $\textnormal{PGL}_2$-invariance of $\widetilde{\omega}_n$ implies that there is a unique algebraic $l$-form $\omega_l$ on $\f{M}_{0,n}$ such that $\omega_l \wedge \nu = \widetilde{\omega}_n$. The dihedral invariance of $\widetilde{\omega}_n$ ensures that $\omega_l$ is invariant under the action of $D_{2n}$ on $\f{M}_{0,n}$.\\
    Explicitly, we may take
    \[
        \nu= \frac{dz_1dz_2dz_3}{(z_1-z_2)(z_2-z_3)(z_3-z_1)},
    \]
    which in turn yields, in simplicial and cubical coordinates,
    \[
        \omega_l=\frac{dt_1\dots dt_l }{ t_2 (t_3 - t_1) (t_4 - t_2) \dots ( t_l-t_{l-2}) (1-t_{l-1} ) } = \frac{dx_1\dots dx_l}{(1-x_1x_2)(1-x_2x_3)\dots (1-x_{l-1}x_l)}.
    \]
    This rightmost expression for $\omega_l$ makes it apparent that $\omega_l$ has neither zeros nor poles along the boundary divisors $\f{M}_{0,n}$ corresponding to $Y_{\lambda_1},\dots, Y_{\lambda_l}\in \F_l$ with $\lambda_i=\{0,\dots, i\}$. The dihedral group acts transitively on the set of divisors at finite distance from $\overline{X}_{n,\delta_0}$ for $\delta_0$ the identity of $\textnormal{Sym}(\{0,\dots, l+1\})$. Since $\omega_l$ is dihedrally invariant, we conclude that $\omega_l$ has neither zeros nor poles at finite distance from $\overline{X}_{n,\delta_0}$. In particular, $\omega_l$ is the unique volume form of $\f{M}_{0,n}$ with this property, up to rational multiples. \\
    In other words, $\omega_l$ is the unique non-vanishing volume form, up to scaling, of the open subscheme of $\overline{\f{M}}_{0,n}$ obtained by removing all boundary divisors that are not at finite distance from $\overline{X}_{n,\delta_0}$. This variety, denoted by $\f{M}_{0,n}^{\delta_0}$, was intensively studied in \cite{Brown-Multiple_zeta_values_and_periods_of_moduli_spaces_of_curves} for the relation of its periods with irrationality proofs for zeta values.  \\
    Strictly speaking, the integral $\xi_l$ is a \emph{cellular integral} in the sense of \cite{Brown-Irrationality_proofs_for_zeta_values_moduli_spaces_and_dinner_parties} only for $l$ even, for it corresponds to the configuration $i \mapsto i+2 \in \text{Sym}_n$. The question of the explicit computation of $\xi_l$ for all $l$ first appeared in \cite[Remark 8.7]{Brown-Multiple_zeta_values_and_periods_of_moduli_spaces_of_curves}.

    \subsection{Some vanishing phenomena}
        
    The goal of this section is to prove the following 
    \begin{proposition}
        \label{prop: weights}
        The integral 
        \[
            \xi_l= \int_{\overline{X}_{n,\delta_0}} \omega_l
        \]
        is a $\Q$-linear combination of multiple zeta values of weight
        \begin{enumerate}
            \item exactly $l$, if $l$ is even;
            \item at most $l-1$, if $l$ is odd.
        \end{enumerate}
    \end{proposition}
    \begin{remark}
        By \cite[Theorem 8.1]{Brown-Irrationality_proofs_for_zeta_values_moduli_spaces_and_dinner_parties}, the integrals $\xi_l$ are $\Q$-linear combinations of multiple zeta values of weight at most $l$.
    \end{remark}
    \noindent
    Although this result will not be strictly necessary for the proof of Theorem~\ref{theo: main}, it is still of independent interest to determine the weights appearing in cellular integrals a priori. We start with a corollary of the study of the polar structure of $\omega_l$ carried out in \cite{Brown-Multiple_zeta_values_and_periods_of_moduli_spaces_of_curves}.
    \begin{lemma}
        \label{lemma: polar structure}
        If $l$ is even, $\omega_l $ has at most simple poles along the boundary divisors of $\overline{\f{M}}_{0,n}$. This also holds if $l$ is odd with the only exception of the divisor associated with $Y_\lambda$ for 
        \[
            \lambda= \{2, 4, \dots, l-1, l+1 \},
        \]
        along which $\omega_l$ has a double pole.
    \end{lemma}
    \begin{proof}
        Let $D$ be the boundary divisor of $\overline{\f{M}}_{0,n}$ associated with the partition $S_1\sqcup S_2=\{z_1,\dots, z_n\}$. By \cite[Proposition 7.5]{Brown-Multiple_zeta_values_and_periods_of_moduli_spaces_of_curves}, we have
        \[
            \textnormal{ord}_D\,  \omega_l = \frac{l-1}{2}-\frac{1}{2} \sum_{i=1}^n \mathbb{I}_D(i,i+2),
        \]
        with indices taken modulo $n$, where
        \[
            \mathbb{I}_D(i,j)= 
            \begin{cases}
                1 & \text{if $\{z_i,z_j\}\subseteq S_1$ or $\{z_i,z_j\}\subseteq S_2$}; \\
                0 & \text{otherwise}.
            \end{cases}
        \]
        Since $\mathbb{I}_D(i,i+2)\le 1$ for all $i=1,\dots, n$, we have 
        \[
            \textnormal{ord}_D\,  \omega_l \ge \frac{l-1}{2}-\frac{n}{2}=-2.
        \]
        Suppose that $\textnormal{ord}_D \, \omega_l=-2$, which implies that $\mathbb{I}_D(i,i+2)=1$ for all $i=1,\dots, n$. In particular, if $i\equiv j \modulo{2}$ for $i,j\in \{1,\dots, n\}$, then $z_i$ and $z_j$ belong to the same set between $S_1$ and $S_2$. Moreover, also $z_n$ and $z_2$ belong to the same set, say $S_1$. \\
        If $l$ is even, so $n$ is odd, this implies that $S_1=\{z_1,\dots, z_n\}$ and $S_2=\emptyset$, against the fact that $\# S_2\ge 2$. It follows that $\textnormal{ord}_D \, \omega_l \ge -1$ for $l$ even. On the other hand, for odd $l$, hence even $n$, we must have $S_1=\{ z_2,z_4, \dots, z_n\}$, while $S_2=\{z_3,z_5,\dots, z_{n-1}\}$. 
    \end{proof}
    \noindent
    The information about the poles of $\omega_l$ provided by Lemma~\ref{lemma: polar structure} is enough to prove Proposition~\ref{prop: weights} in the case of even $l$.
    \begin{lemma}
        \label{lemma: weights for even l}
        If $l$ is even, then $\xi_l$ is a $\Q$-linear combination of multiple zeta values of weight exactly $l$.
    \end{lemma}
    \begin{proof}
        We follow the strategy of \cite[Proposition 3.12]{Dupont-Odd_zeta_motive_and_linear_forms_in_odd_zeta_values}. Since $\omega_l$ is closed, by Lemma~\ref{lemma: polar structure} both $\omega_l$ and $d\omega_l=0$ have at most simple poles along the boundary divisors of $\overline{\f{M}}_{0,n}$. This means that $\omega_l$ has logarithmic singularities along the boundary divisors. \\
        Write for brevity $Y=\overline{\f{M}}_{0,n}$, let $D$ be the polar divisor of $\omega_l$ in $Y$ and let $Z$ be the union of the boundary divisors at finite distance from $\overline{X}_{n,\delta_0}$. Note that $D\cup Z$ has simple normal crossings. The integral $\xi_l$ is a period of the mixed Tate Hodge structure of the cohomology group $H^l_{\textnormal{dR}}(Y\setminus D, Z\setminus(D\cap Z))$. Since the highest term in the Hodge filtration reads
        \[
            F^lH^l_{\textnormal{dR}}(Y\setminus D, Z\setminus(D\cap Z))= \textnormal{Im}(H^0(\Omega^{l}_Y(\log D))\to  H^l_{\textnormal{dR}}(Y\setminus D, Z\setminus(D\cap Z))),
        \]
        the cohomology classes of pure weight $2l$ are precisely those with logarithmic singularities along $D$, hence $\omega_l$ is one of these. This implies that $\xi_l$ is a $\Q$-linear combination of multiple zeta values of weight exactly $l$.
    \end{proof}
    \begin{remark}
        The differential form $\omega_l$ for $l$ even is a \emph{cell-form} on $\f{M}_{0,n}$. A study of the $\Q$-algebra generated by these forms has been conducted in \cite{Brown-Carr-Schneps-The_algebra_of_cell-zeta_values}.
    \end{remark}
    \noindent
    The drop of the weight in the case of odd $l$ is due to the existence of an algebraic primitive for $\omega_l$. This is not a consequence of Lemma~\ref{lemma: polar structure}, but the latter can help in finding a primitive, as we will explain in more detail in the next section.\\
    For brevity, let us set 
    \[
        f=(1-x_1x_2)(1-x_2x_3)\dots (1-x_{l-1}x_l).
    \]
    The following lemma completes the proof of Proposition~\ref{prop: weights}
    \begin{lemma}
    \label{lemma: weight for odd l}
    Consider the following $(l-1)$-form on $\f{M}_{0,n}$:
        \[
            \alpha_l=\sum_{i=1}^l \frac{x_i}{f} \, dx_1\dots dx_{i-1}dx_{i+1}\dots dx_l.
        \]
        Then we have 
        \begin{enumerate}
            \item $d \alpha_l=\omega_l$ if $l$ is odd;
            \item $d\alpha_l=0$ if $l$ is even.
        \end{enumerate}
        In particular, if $l$ is odd, then $\xi_l$ is a $\Q$-linear combination of multiple zeta values of weight at most $l-1$.
    \end{lemma}
    \begin{proof}
        The differential of $\alpha_l$ is
        \begin{align*}
            d\alpha_l &= \left(\sum_{i=1}^l (-1)^{i+1} \frac{\partial}{\partial x_i} \left(\frac{x_i}{f}\right) \right) \, dx_1\dots dx_l = \left(\sum_{i=1}^l (-1)^{i+1} \left(\frac{1}{f}-\frac{x_i}{f^2}\frac{\partial f}{\partial x_i} \right) \right) \, dx_1\dots dx_l \\
                    &= \left(\sum_{i=1}^l \frac{(-1)^{i+1}}{f}  + \frac{1}{f^2} \sum_{i=1}^l x_i \frac{\partial f}{\partial x_i} \right) \, dx_1\dots dx_l.
        \end{align*} 
        The first sum equals $f^{-1}$ if $l$ is odd, while it vanishes otherwise. For the second summand, note that for $i\neq 1, l$, we have 
        \[
            \frac{\partial f}{\partial x_i} = \left(- \frac{x_{i-1}}{(1-x_{i-1}x_{i})}-\frac{x_{i+1}}{(1-x_ix_{i+1})} \right) f,
        \]
        while the derivatives with respect to $x_1$ and $x_l$ read
        \[
            \frac{\partial f}{\partial x_i}=-\frac{x_2}{1-x_1x_2}f, \qquad \frac{\partial f}{\partial x_l}=-\frac{x_l}{1-x_{l-1}x_l}f.
        \]
        Hence, by a telescopic summation, we have
        \[
            \sum_{i=1}^l (-1)^i x_i \frac{\partial f}{\partial x_i}=0,
        \]
        which proves the claim about the differential of $\alpha_l$. \\
        Suppose that $l$ is odd, so $\alpha_l$ is an algebraic primitive of $\omega_l$. Since $\omega_l$ has no poles along the boundary divisors at finite distance from $\overline{X}_{n,\delta_0}$, there is an algebraic primitive $\widetilde{\alpha}_l$ of $\omega_l$ that also has no poles along along these divisors. To obtain $\widetilde{\alpha}_l$, one may apply to $\alpha_l$ the regularization procedure explained in \cite[Proposition 8.1]{Brown-Multiple_zeta_values_and_periods_of_moduli_spaces_of_curves}. The so-constructed primitive remains algebraic because this procedure does not increase the weight. \\
        The absence of poles on the boundary of $\overline{X}_{n,\delta_0}$ ensures that $\widetilde{\alpha}_l$ restricts continuously to the involved boundary divisors. By Stokes' theorem,
        \[
            \int_{\overline{X}_{n\delta_0}}\omega_l=\int_{\partial\overline{X}_{n\delta_0}} \widetilde{\alpha}_l.
        \]
        The boundary divisors of $\overline{\f{M}}_{0,n}$ are isomorphic to products of the form $\overline{\f{M}}_{0,n_1}\times \overline{\f{M}}_{0,n_2}$ with $n_1+n_2=n-1$, whose periods are linear combinations of multiple zeta values of weight at most $l-1$. The claim follows.
    \end{proof}

    \subsection{A case of study in three dimensions}
        
    If $l$ is odd, we have seen that $\omega_l$ is exact as an algebraic form. However, in general, finding an explicit primitive can be a hard task. We will explain in this section how the observation about the poles of $\omega_l$ exposed in Lemma~\ref{lemma: polar structure} can help in this context. To simplify the computations, we will focus on the case $l=3$. \\
    Recall that, by Lemma~\ref{lemma: polar structure}, $\omega_3$ has at most simple poles, with the only exception of a double pole along the divisor $D_X$ associated with the irreducible element $X=\langle y_2\rangle \in \mathcal{F}_3$. Let us choose a local chart on $\overline{\f{M}}_{0,n}$ for which this double pole becomes visible. Consider the maximal $\F_3$-nested set $\mathcal{S}=\{ \langle y_2\rangle, \langle y_2,y_4\rangle, \langle y_1, y_2, y_4\rangle, V^*  \}$ with adapted basis
    \[
        \langle y_2\rangle \mapsto y_2, \quad \langle y_2,y_4\rangle \mapsto y_2-y_4, \quad \langle y_1, y_2, y_4\rangle \mapsto y_1, \quad V^*\mapsto y_3.
    \]
    The closure of the subscheme $\{u_2=0\}\subseteq U_\S$ in $\overline{\f{M}}_{0,n}$ coincides with $D_X$. The induced standard coordinates $u_1,u_2,u_3$ on the local chart $U_\mathcal{S}$ are given by
    \[
    u_1=\frac{y_1}{y_3}, \quad u_2=\frac{y_2}{y_2-y_4}, \quad u_3=\frac{y_2-y_4}{y_1},
    \]
    with inverse
    \[
    y_1= u_1, \quad, y_2=u_1u_2u_3, \quad y_4=(u_2-1)u_1u_3.
    \]
    In these coordinates, we see that 
    \[
        \omega_3= \frac{1}{u_3^2}\cdot \frac{1}{u_1(1-u_1)u_2(1-u_2)}\, du_1 du_2 du_3.
    \]
    At this point, it is straightforward to find an algebraic primitive for $\omega_3$, namely
    \[
        \alpha_3=-\frac{1}{u_3}\cdot \frac{1}{u_1(1-u_1)u_2(1-u_2)} \, du_1 du_2.
    \]
    Going back to simplicial and cubical coordinates, we compute
    \begin{align*}
        \alpha_3&=\frac{t_3}{t_2(t_3-t_1)(1-t_2)}\, dt_1dt_2 + \frac{t_1}{t_2(t_3-t_1)(1-t_2)} \, dt_2 dt_3 \\
        &= \frac{x_1}{f}\, dx_2dx_3 + \frac{x_2}{f}\, dx_1dx_3 + \frac{x_3}{f}\, dx_1dx_2.
    \end{align*}
    From this expression it is not too difficult to guess the shape of the primitive $\alpha_l$ found in the previous section. 

    \begin{remark}
        If we turn the double pole along $D_X$ into a simple one, by the same argument as in Lemma~\ref{lemma: weights for even l} we obtain multiple zeta values only of weight $3$. To achieve this, it is not enough to multiply $\omega_3$ by $u_3$, because
        \[
            \frac{u_3}{y_4}=\frac{t_2-1}{t_1}=\frac{x_2x_3-1}{x_1x_2x_3}
        \]
        would introduce poles along the boundary of the integration domain. We could instead multiply by
        \[
            u_1u_2u_3=\frac{y_2}{y_3}=\frac{t_2}{t_3}=x_2.
        \]
        In higher dimensions, for the first few values of $l$ odd, we have checked by means of the program \texttt{HyperInt} \cite{Panzer-Algortihms_for_the_symbolic_integration_of_hyperlogarithms_with_applications_to_Feynman_integrals} that the integral
        \[
            \int_{[0,1]^l} x_2\omega_l
        \]
        is a linear combination of multiple zeta values of weight exactly $l$. By symmetry, this also applies to $x_i\omega_l$ for all $2\le i\le l-1$ even.
    \end{remark}
    $ $\\
    
    Recall that Theorem~\ref{theo: main} predicts the relation $\xi_3=2\xi_2$. We explain how to derive this inequality using the primitive $\alpha_l$ in the spirit of the period conjecture. In the next section, we set up the argument for the general formula for odd $l$. \\
    We follow the general strategy of Lemma~\ref{lemma: weight for odd l}. First, we need to regularize the primitive $\alpha_3$ along the boundary of $\Delta=\overline{X}_{n,\delta_0}$. The irreducible divisors at finite distance correspond to the following subset of $\F_3$:
    
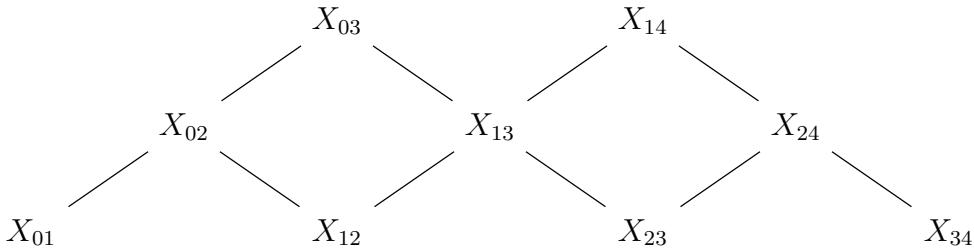
\[
    \begin{tikzcd}
        & & X_{03} & & X_{14} & & \\
         &X_{02} \ar[ur,dash] & & X_{13} \ar[ur,dash] \ar[ul, dash] & & X_{24} \ar[ul, dash]&  \\
    X_{01} \ar[ur,dash] & & X_{12} \ar[ur,dash] \ar[ul, dash] & & X_{23} \ar[ur,dash]\ar[ul, dash] & & X_{34}\ar[ul, dash]
    \end{tikzcd}
    \]
    Here, for $0<i< j$ we have written for short $X_{ij}=\langle y_k-y_h \mid i\le k,h \le j \rangle$, while $X_{0j}=\langle y_k\mid 1\le k \le j\rangle$. To describe the poles of $\alpha_3$ along the divisors $D_{X_{ij}}$, we choose three maximal $\F_3$-nested sets which make up a partition of the $X_{ij}$'s:
    \begin{enumerate}
        \item $\S_1=\{X_{01}, X_{02}, X_{34}, V^* \} $;
        \item $\S_2=\{X_{12}, X_{13}, X_{03}, V^* \} $;
        \item $\S_3=\{X_{23}, X_{24}, X_{14}, V^* \} $.
    \end{enumerate}
    As adapted bases, we choose those induced by the assignment $X_{0j}\mapsto y_j $ and $X_{ij}\mapsto y_i-y_j$ for $i\neq 0$. The standard coordinate in the chart $U_{S_k}$ corresponding to $X_{ij}$ will be denoted by $v_{ij}$.
    By expressing $\alpha_3$ in each of the charts $U_{\S_k}$, it turns out that it has simple poles along the divisors corresponding to $X_{13}$ and $X_{14}$.\\
    We may regularize $\alpha_3$ by removing, one after the other, the singular part of its Laurent expansion at each divisor. We will explain this in more detail in the next section; for the moment, this method leads to the primitive
    \begin{align*}
        \widetilde{\alpha}_3 =& \, \frac{t_1t_2-t_3t_1-t_2t_3+t_3}{t_2(t_3-t_1)(1-t_2)(1-t_1)}\, dt_1dt_2 \\
                & \, + \frac{t-1}{t_2(t_3-t_1)(1-t_2)}\, dt_2 dt_3- \frac{1}{(1-t_1)(t_3-t_1)}dt_1dt_3. 
    \end{align*}
    By Stokes' theorem, we then have
    \[
    \xi_3=\int_{\Delta} \omega_3 = \int_{\partial\Delta } \widetilde{\alpha}_3.
    \]
    We want to compute the restriction of $\widetilde{\alpha}_3$ to each  boundary divisor at finite distance with $\Delta=\overline{X}_{n,\delta_0}$. Writing  $\widetilde{\alpha}_3$ in the three local charts that we have isolated, it turns out that the boundary divisors where $\widetilde{\alpha}_3$ does not vanish are given by $X_{34}$, where $\widetilde{\alpha}_3$ restricts to
        \[
            \widetilde{\alpha}_3\vert_{v_{34}=0}=\frac{1}{1-v_{01}v_{02}} \, dv_{01} dv_{02},
        \]
        and $X_{14}$, where $\widetilde{\alpha}_3$ restricts to
        \[
            \widetilde{\alpha}_3\vert_{v_{14}=0}=\frac{1}{1-v_{24}+v_{23} v_{24}} \,dv_{23} dv_{24}.
        \]
    Thus, we have
    \[
    \xi_3= \int_{\partial\Delta } \widetilde{\alpha}_3= \int_{\partial\Delta \cap X_{34}}\frac{1}{1-v_{01}v_{02}} \, dv_{01} dv_{02} + \int_{\partial \Delta \cap X_{14}}  \frac{1}{1-v_{24}+v_{23} v_{24}} \,dv_{23} dv_{24}.
    \]
    It is clear that $\partial \Delta\cap X_{34}$ is given by $[0,1]^2$, simply by looking at the coordinates on $U_{\S_1}$. On the other hand, the restriction of the arrangement $\A_3$ to $X_{14}$ is isomorphic to $\A_2$ by setting $y_1'=y_1-y_2$, $y_2'=y_1-y_3$ and $y_3'=y_1-y_4$. In the affine coordinates $t_i'=\frac{y_i'}{y_3'}$ we have
    \[
        v_{23}=\frac{t_2'-t_1'}{1-t_1'}, \qquad v_{24}=1-t_1'.
    \]
    It follows that
    \[
        \frac{1}{1-v_{24}+v_{23} v_{24}} \,dv_{23} dv_{24} = \frac{1}{t_2'(1-t_1')}\, dt_1'dt_2'.
    \]
    It is readily checked that the integral takes place over the standard simplex, hence we conclude that also the second integral equals $\xi_2$. \\
    Thus, the formula $\xi_3=2\xi_2$ can be explained by means of Stokes' theorem and the restrictions of $\widetilde{\alpha}_3$ to the boundary divisors, which are either $0$ or $\omega_2$. In general, we would expect that the formula
    \[
        \xi_{2m+1}=\sum_{k=0}^m \xi_{2k}\xi_{2m-2k}
    \]
    can be explained in a similar way, by finding a regularized primitive $\widetilde{\alpha}_l$ of $\omega_l$ whose restriction to the divisors in the boundary is either $0$ or $\omega_{2k}\omega_{l-1-2k}$. This product would be justified by the decomposition of boundary divisors as products of moduli spaces $\f{M}_{0,n'}$ of smaller dimension.

    \subsection{Polar structure of the primitive}
        
    For this section, we fix $l$ odd, hence $n=l+3$ is even. Giving an explanation in terms of the period conjecture of the formula expressing the odd-dimensional integrals with respect to the even-dimensional ones appears quite intricate. Despite having a simple shape, the primitive $\alpha_l$ is not particularly well-behaved under the action of the automorphisms of $\overline{\f{M}}_{0,n}$. This makes the regularization process very difficult to handle. \\
    To shed more light on the problem, we expose the polar structure of $\alpha_l$ in this section. This is the only place where we make use of the full power of the moduli description of $\overline{\f{M}}_{0,n}$ and the symmetries coming with it.\\
    Recall that the boundary divisors at finite distance from $\overline{X}_{n,\delta_0}$ are in bijection with the partitions of the set $\{z_1,\dots, z_n\}$ into two elements of cardinality at least $2$. Upon identifying this set with the edges of an $n$-gon along the standard order induced by $\delta_0$, these partitions correspond to the diagonals of the $n$-gon. Let us denote by $D_{ij}$ the divisor associated with the diagonal $\{i,j\}$ connecting the vertices between the edges $i,i+1$ and $j,j+1$. When referring to these diagonals, we will always take indices modulo $n$. \\
    The goal of this section is to prove the following:
    \begin{proposition}
        \label{prop: poles of the primitive}
        Among the boundary divisors at finite distance from $\overline{X}_{n,\delta_0}$, the differential form $\alpha_l$ has a pole precisely at the divisors
        \begin{enumerate}
            \item $D_{13},\dots, D_{1n-2}$;
            \item $D_{1-h,4-h},\dots D_{1-h,1}$ for $h=3,\dots, n-3$;
            \item $D_{1-h,3},\dots, D_{1-h,n-2-h}$ for $h=1,\dots, n-5$.
        \end{enumerate}
        These poles are all simple.
    \end{proposition}
    \noindent
    Recall that $D_{i,j}\cap D_{i',j'}\neq \emptyset$ if and only if the two associated irreducible subsets of $L(\A_l)$ are $\F_l$-nested; in terms of the $n$-gon, this means that the diagonals $\{i,j\}$ and $\{i',j'\}$ do not cross. A maximal $\F_l$-nested set corresponds therefore to a triangulation of the $n$-gon by its diagonals.\\
    Let $D$ denote the union of the boundary divisors that are not at finite distance from $\overline{X}_{n,\delta_0}$. We describe dihedral coordinates on $\f{M}_{0,n}^{\delta_0}=\overline{\f{M}}_{0,n}\setminus D$ following \cite{Brown-Multiple_zeta_values_and_periods_of_moduli_spaces_of_curves}. Consider symbols $u_{ij}$ for $i,j=1,\dots, n$, $i\neq j-1,j,j+1$ subject to the relation $u_{ij}=u_{ji}$. Given a set of diagonals $A$ of the $n$-gon, we define $u_A$ as the product of all $u_{ij}$ for which the diagonal ${i,j}$ belongs to $A$.
    \begin{lemma}[{\cite[Lemma 2.30]{Brown-Multiple_zeta_values_and_periods_of_moduli_spaces_of_curves}}]
        There is an isomorphism
        \[
            \f{M}_{0,n}^{\delta_0}\cong\textnormal{Spec}\, \Q[u_{ij}]/I^\chi_{n,\delta_0},
        \]
        where the ideal $I^\chi_{n,\delta_0}$ is generated by elements of the form $u_{ij}-u_{ji}$ and $u_A+u_B=1$ for all sets of diagonals $A,B$ that cross completely (see \cite[Section 2.2]{Brown-Multiple_zeta_values_and_periods_of_moduli_spaces_of_curves}). The closure of the subscheme $\{u_{ij}=0\}$ in $\overline{\f{M}}_{0,n}$ coincides with $D_{ij}$.
    \end{lemma}
    \noindent 
    We exploit this description to provide some local charts on $\f{M}_{0,n}^{\delta_0}$ that allow us to expose the polar structure of $\alpha_l$. For all fixed $i=1,\dots, n$ the set $A_i$ of diagonals $\{i,j\}$ with $j\neq i-1,i,i+1$ provides a full triangulation of the $n$-gon. By inverting all coordinates $u_{i'j'}$ for diagonals $\{i',j'\}$ not belonging to $A_i$ one recovers the local chart $U_{\S_i}$ for the maximal $\F_l$-nested set $\S_i$ corresponding to $A_i$. In particular, we have $u_{2j}=x_j$ for $j=4,\dots, n$. \\
    The relations provided by $I^\chi_{n,\delta_0}$ give
    \[
        \begin{dcases}
            u_{2n}=1-u_{13}\dots u_{1n-1} \\
            u_{2j}= \frac{1-u_{13}\dots u_{1j-1}}{1-u_{13}\dots u_{1j}} & (j=4,\dots, n-1).
        \end{dcases}
    \]
    Let $\rho$ be the automorphism of $\overline{\f{M}}_{0,n}$ associated with the dihedral symmetry $i \mapsto i+1$, which acts as a rotation of the $n$-gon. This induces an automorphism of $\f{M}_{0,n}^{\delta_0}$ that sends $u_{ij}$ to $u_{i-1,j-1}$. \\
    To study the poles of $\alpha_l$, we express $\alpha_l$ with respect to $u_{13},\dots, u_{1n-1}$ via the above change of coordinates. Then, we apply $\rho$ and once again the same change of coordinates. Repeating this $n$-times, we go through a set of local charts in which every boundary divisor becomes visible at least once. \\

    To make notation lighter, let us set for all $k=3,\dots, n-1$ and $j\ge k$
    \[
        du_{13}\overset{k}{\dots} du_{1j}=du_{13}\dots du_{1k-1} du_{1k+1} \dots du_{1j}
    \]
    If $j=n-1$, we write $\lambda_k=du_{13}\overset{k}{\dots} du_{1n-1}$. \\
    It is straightforward to check that for $j=4,\dots, n-1$
    \begin{align*}
        du_{24}\dots du_{2j} =& \frac{\prod_{h=3}^{j-2}u_{1h}^{j-h-1}}{(1-u_{13}u_{14})\dots (1-u_{13}\dots u_{1j-1})(1-u_{13}\dots u_{1j})^2}\times \\
        & \times \bigg( \sum_{k=3}^{j-1} (-1)^k u_{1k}\dots u_{1j-1}(u_{1k}-1) du_{13}\overset{k}{\dots} du_{1j} +\\
        & \qquad  + (-1)^j (u_{1j}-1)du_{13}\dots du_{1j-1} \bigg).
    \end{align*}
    Moreover, for $j=4,\dots, n-1$
    \begin{align*}
        du_{2j+1}\dots du_{2n} 
        &= \frac{(-1)^{n-j}u_{13}^{n-j}\dots u_{1j}^{n-j}\prod_{h=j+1}^{n-2} u_{1h}^{n-h-1}}{(1-u_{13}\dots u_{1j+1})\dots (1-u_{13}\dots u_{1_{n-1}})} \left( \sum_{k=3}^j \frac{du_{1k}du_{1j+1}\dots du_{1n-1}}{u_{1k}} \right).
    \end{align*}
    Recall that in the last section we have defined
    \[
        f= \prod_{k=4}^{n-1}(1-u_{2k}u_{2k+1})=\prod_{k=3}^{n-2} u_{1k}^{n-k-1} \prod_{j=4}^{n-2} \frac{1-u_{1j}u_{1j+1}}{1-u_{13}\dots u_{1j+1}}.
    \]
    For $j=4,\dots, n-1$, we then have
    \begin{align*}
        \frac{u_{2j}}{f}du_{24}\overset{j}{\dots}du_{2n} &=  \frac{(-1)^{n-1}}{\rho^*f }\frac{(u_{13}\dots u_{1j-1}-1)u_{1j}}{u_{13}\dots u_{1j-1}}du_{13}\overset{j}{\dots} du_{1 n-1} +&\\
        &\qquad +\sum_{k=3}^{j-1} \frac{(-1)^{n-j+k}}{f((u_{1k})_k)} \frac{(u_{1k}-1)}{u_{13}\dots u_{1k-1}} du_{13}\overset{k}{\dots} du_{1n-1}&
    \end{align*}
    This formula also works for $j=n$, provided that we formally set $u_{1n}=0$. Since $n$ is even, it follows that
    \begin{align*}
        \alpha_l =& \frac{1}{\rho^*f} \bigg( \sum_{\substack{k=3 \\ \text{$k$ odd}}}^{n-1}\frac{1-u_{13}\dots u_{1k}}{u_{13}\dots u_{1k-1}} \lambda_k + \sum_{\substack{k=3 \\ \text{$k$ even}}}^{n-1}\frac{(1-u_{13}\dots u_{1k-1})u_{1k}}{u_{13}\dots u_{1k-1}} \lambda_k\bigg).
    \end{align*}
    Let us now write for short $\gamma_l=\rho^* \alpha_l$ and $\gamma_l^{(h)}=(\rho^{-h})^*\gamma_l$. Proposition~\ref{prop: poles of the primitive} follows at once from the following description of $\gamma^{(h)}_l$ using that $\alpha_l=(\rho^h)^*\gamma^{(h+1)}$.
    \begin{lemma}
    \label{lemma: primitive in all coordinates}
        We have
        \[
            \gamma_l^{(h)}= \sum_{k=3}^{n-1} \frac{p_{h,k}}{\rho^*f} \lambda_k,
        \]
        where the coefficients $p_{h,k}$ are defined as follows.
        For $h=1$ we have 
        \[
            p_{1,k}= 
            \begin{dcases}
                \frac{1-u_{13}\dots u_{1k}}{u_{13}\dots u_{1k-1}} & \text{if $k$ is odd;} \\
                \frac{u_{1k}(1-u_{13}\dots u_{1k-1})}{u_{13}\dots u_{1k-1}} & \text{if $k$ is even.}
            \end{dcases}
        \]
        For $h=2,\dots n-2$ and $k$ odd, we have:
        \[
            p_{h,k}= 
            \begin{dcases}
                \frac{1+u_{1k}\dots u_{1h}(u_{1k+1}\dots u_{1h+1}-u_{1h+1}-1)}{u_{1k+1}\dots u_{1h}(1-u_{1h+1})} & \text{if $k\le h$;} \\
                u_{1k} & \text{if $k=h+1$;} \\
                 \frac{1+u_{1h+2}\dots u_{1k}(u_{1h+1}\dots u_{1k-1}-u_{1h+1}-1)}{u_{1h+2}\dots u_{1k-1}(1-u_{1h+1})} & \text{ if $k\ge h+2$.}
            \end{dcases}
        \]
        For $h=2,\dots n-2$ and $k$ even, we have:
        \[
            p_{h,k}= 
            \begin{dcases}
                \frac{u_{1k}(1+u_{1k+1}\dots u_{1h}(u_{1k}\dots u_{1h+1}-u_{1h+1}-1))}{u_{1k+1}\dots u_{1h}(1-u_{1h+1})} & \text{if $k\le h$;} \\
                -u_{1k} & \text{if $k=h+1$;} \\
                 \frac{u_{1k}( 1+u_{1h+2}\dots u_{1k-1}(u_{1h+1}\dots u_{1k}-u_{1h+1}-1))}{u_{1h+2}\dots u_{1k-1}(1-u_{1h+1})} & \text{ if $k\ge h+2$.}
            \end{dcases}
        \]
        For $h=n-1$ we have
        \[
            p_{n-1,k}=
            \begin{dcases}
                \frac{1-u_{1k}\dots u_{1n-1}}{u_{1k+1}\dots u_{1n-1}} & \text{if $k$ is odd;} \\
                \frac{u_{1k}(1-u_{1k+1}\dots u_{1n-1})}{u_{1k+1}\dots u_{1n-1}} & \text{if $k$ is even.}
            \end{dcases}
        \]
    \end{lemma}
    \begin{proof}
        Define $\epsilon_k=1$ for $k$ even, $\epsilon_k=0$ for $k$ odd. Consider the differential forms
        \begin{equation*}
            \phi_l=  \frac{1}{\rho^*f} \sum_{k=3}^{n-1}\frac{u_{1k}^{\epsilon_k}}{u_{13}\dots u_{1k-1}} \lambda_k, \qquad
            \chi_l =\frac{1}{\rho^*f} 
        \sum_{ k=3}^{n-1}
        u_{13}\dots u_{1k}^{1+\epsilon_k} \lambda_k.
        \end{equation*}
        The computations carried out before let one easily check that 
        \begin{align*}
                (\rho^{-1})^*\gamma_l &=-\gamma_l+\phi_l, \\
                (\rho^{-1})^*\phi_l &=-\frac{2}{1-u_{13}}\gamma_l+\frac{1}{1-u_{13}}\phi_l+ \frac{1}{1-u_{13}}\chi_l, \\
                (\rho^{-1})^* \chi_l &= (1-u_{13}) \phi_l. 
        \end{align*}
        In particular, we may write $\gamma_l^{(h)}= a_h \gamma_l +b_h \phi_l + c_h\chi_l$ for suitable regular functions $a_h$, $b_h$ and $c_h$. To find the latter, we observe that from the equations above
        \begin{align*}
            a_h &= -(\rho^{-1})^*a_{h-1}-\frac{2}{1-u_{13}}\,(\rho^{-1})^*b_{h-1}, \\
            b_h &= (\rho^{-1})^*a_{h-1}+\frac{1}{1-u_{13}}\,(\rho^{-1})^*b_{h-1}+(1-u_{13})\, (\rho^{-1})^*c_{h-1}, \\
            c_h &= \frac{1}{1-u_{13}}\, (\rho^{-1})^*b_{h-1}.
        \end{align*}
        It is now straightforward to check by induction that 
        \begin{align*}
            a_0&=1, & b_0&=0, & c_0&=0; \\
            a_1&=-1, & b_1&=1, & c_1&=0; \\
            a_h&=\frac{u_{1h+1}+1}{u_{1h+1}-1}, & b_h&=\frac{u_{13}\dots u_{1h+1}}{1-u_{1h+1}}, & c_h&=\frac{1}{u_{13}\dots u_{1h}(1-u_{1h+1})}; \\
            a_{n-1}&=-1, & b_{n-1}&=0, & c_{n-1}&=\frac{1}{u_{13}\dots u_{1n-1}}; \\
            a_n&=1, & b_n&=0, & c_n&=0; 
        \end{align*}
        where $h=2,\dots, n-2$. The claim follows at once.
    \end{proof}
    \noindent
    In order to apply Stokes' theorem and turn $\xi_l$ for odd $l$ into a sum of periods of $\overline{\f{M}}_{0,n'}$ with $n'<n$, we need to regularize $\alpha_l$ by removing its poles at finite distance form $\overline{X}_{n,\delta_0}$. We follow the procedure of \cite[Proposition 8.1]{Brown-Multiple_zeta_values_and_periods_of_moduli_spaces_of_curves} and \cite{Panzer-Algortihms_for_the_symbolic_integration_of_hyperlogarithms_with_applications_to_Feynman_integrals}.\\
    Suppose we want to remove a pole of $\alpha_l$ along a boundary divisor $D_{ij}$ at finite distance form $\overline{X}_{n,\delta_0}$. We first express $\alpha_l$ in Laurent series with respect to the coordinate $u_{ij}$. Then, we consider the differential form $\widetilde{\alpha}_l$ obtained from $\alpha_l$ by removing the singular part of this series expansion. One can check that $d\widetilde{\alpha}_l=\omega_l$ exploiting the smoothness of $\omega_l$ along $D_{ij}$. Moreover, this procedure does not produce further poles along other boundary divisors. Thus, we may repeat it in turn for all divisors along which $\alpha_l$ has a pole until we obtain a primitive that can be restricted to $\partial \overline{X}_{n,\delta_0}$. \\
    Unfortunately, the polar structure of $\alpha_l$ is highly asymmetric, so the regularization procedure becomes unmanageable very quickly. To exemplify this, we regularize $\alpha_l$ with respect to $u_{13},\dots, u_{1n-2}$. Following the notation of \cite{Brown-Multiple_zeta_values_and_periods_of_moduli_spaces_of_curves}, the only regularizations that we need are
    \begin{align*}
        \textnormal{Reg}\left(\frac{1}{\rho^*f \cdot  u_{13}} , D_{13}\right) &= \frac{1}{\rho^*f \cdot u_{13}}- \frac{1-u_{13}u_{14}}{\rho^*f \cdot u_{13}}= \frac{u_{14}}{\rho^*f }, \\
        \textnormal{Reg}\left(\frac{1}{\rho^*f \cdot u_{1q}} , D_{1q}\right) &= \frac{1}{\rho^*f \cdot u_{1q}}- \frac{(1-u_{1q-1}u_{1q})(1-u_{1q}u_{1q+1})}{\rho^*f \cdot u_{1q}}\\
        &= \frac{u_{1q-1}}{\rho^*f }+\frac{u_{1q+1}(1-u_{1q-1}u_{1q})}{\rho^*f},
    \end{align*}
    for $q=4,\dots, n-2$. Define $X_4=0$, $X_5=1$ and recursively for $q\ge 6$
    \[
        X_q=u_{1q-2}X_{q-1}+(1-u_{1q-3}u_{1q-2})X_{q-2},
    \]
    which is a polynomial in $u_{13},\dots, u_{1q-2}$. By iteratively applying the above regularizations to $\alpha_l$, one obtains a differential form $\widetilde{\alpha}_l$ of the form
    \begin{align*}
        \widetilde{\alpha}_l=& -\rho^*\alpha_l+  \frac{1}{\rho^*f}\biggl( \lambda_3+u_{14}^2\lambda_4+\lambda_5+u_{16}(u_{16}+u_{14}-u_{14}u_{15}u_{16})\lambda_6 \\
        &+ \sum_{k=7}^{n-1} u_{1k}^{\epsilon_{k}} ( X_{k-2}(1-u_{1k-3}u_{1k-2})+X_{k-1}( u_{1k}+u_{1k-2}-u_{1k-2}u_{1k-1}u_{1k}))\lambda_k\biggr).
    \end{align*}
    After this, the necessary changes of coordinates and regularizations become quite unfeasible to handle. 
    \begin{remark}
        It would be interesting to have access to an equivariant regularization procedure. Note that the differential form
        \[
            \frac{1}{n}\sum_{h=1}^n \gamma_l^{(h)} =\frac{1}{n} \sum_{h=1}^n (\rho^h)^*\alpha_l
        \]
        is invariant under the rotation $\rho$. Finding a regularization of this differential form with respect to a single boundary divisor that is again $\rho$-invariant would allow us to remove several poles at once.
    \end{remark}

    \section{Main result}
    \subsection{Preliminary definitions}
        
    The rest of the article is devoted to the proof of Theorem~\ref{theo: main}. The method is an inductive application of the strategy first proposed in \cite{Brown-Multiple_zeta_values_and_periods_of_moduli_spaces_of_curves} and later turned into an effective algorithm in \cite{Panzer-Algortihms_for_the_symbolic_integration_of_hyperlogarithms_with_applications_to_Feynman_integrals}. For a given $l\ge 2$, all the computations can be checked with the program \texttt{HyperInt} by Eric Panzer \cite{Panzer-Algortihms_for_the_symbolic_integration_of_hyperlogarithms_with_applications_to_Feynman_integrals}. \\ 
    We first recall some notions about hyperlogarithms and polylogarithms. Consider a one-dimensional hyperplane arrangement over $\C$, corresponding to $\mathbb{P}^1$ with the points $\sigma_0,\sigma_1,\dots, \sigma_n\in k$ and $\infty$ removed. We may assume that $\sigma_0=0$ and $\sigma_1=1$.
    Let $\omega_i=d\log (z-\sigma_i)$ and let $B$ be the set of non-commutative polynomials over $\Q$ in the $\omega_i$'s; we will use the bar notation $[\omega_{i_s}\vert\dots \vert\omega_{i_{0}}]$ for standard monomials in $B$, which are usually referred to as \emph{words}. We endow $B$ with the \emph{shuffle product} $\sh$, defined as follows. Given monomials $v,w\in B$ and $i,j\in\{0,\dots, n\}$, we set inductively
    \[
        [\omega_i\vert v] \,\sh\,[\omega_j\vert w]= \left[\omega_i\vert (v \,\sh\, [\omega_j\vert  w])\right]+\left[\omega_j\vert ([\omega_i \vert v] \,\sh\, w)\right],
    \]
    together with $ w\,\sh\, 1=1 \,\sh\, w = w$ and $1 \,\sh\, 1 =1$. This definition extends to the whole $B$ by $\Q$-linearity. The resulting commutative $\Q$-algebra $B$ is the \emph{reduced bar complex} of $\mathbb{P}^1\setminus \{0,\sigma_1,\dots, \sigma_n,\infty\}$.  \\
    For every monomial $w=[\omega_{i_s}\vert \dots\vert \omega_{i_0}]\in B $, its associated \emph{hyperlogarithm} $L_w(z)$ is defined inductively as follows. If $w=[\omega_0^m]$ is the monomial given by repeating $\omega_0$ $m$-times, then 
    \[
        L_w(z)=\frac{1}{m!}\log^m z.
    \]
    Otherwise, if $i_0\neq 0$ we define inductively the multi-valued function
    \[
        L_w(z)=\int_0^z\frac{L_{[\omega_{i_{s-1}}\vert \dots \vert \omega_{i_0}]}(t)}{\sigma_{i_s}-t}\,dt.
    \]
    To the empty word we associate by definition $L_\varnothing (z)=1$. This represents $L_w(z)$ as the iterated integral
    \[
        L_w(z)=\int_0^z \frac{1}{\sigma_{i_s}-t_s}\int_0^{t_s}\frac{1}{\sigma_{i_{s-1}}-t_{s-1}}\dots \int_0^{t_{1}}\frac{1}{\sigma_{i_0}-t_0}dt_0\dots dt_s,
    \]
    which is convergent because $\sigma_{i_0}\neq 0$. To extend the map $w\mapsto L_w(z)$ to the whole $B$ as a $\Q$-algebra homomorphism, it suffices to observe that every standard monomial $w\in B$ can be written as 
    \[
        w=\sum_{v}\sum_{m\ge0} a_{v,m} v \,\sh\, [\omega_0^m]
    \]
    for some uniquely determined $a_{v,n}\in \Q$, where the first sum runs over all monomials $v$ which do not end with $\omega_0$.\\
    We may describe these functions locally around zero in terms of multiple polylogarithms. Take $r\ge 1$ and let $n_1,\dots, n_r\in \Z$ be positive integers. The multiple polylogarithm associated with $(n_1,\dots, n_r)$ is defined in a neighborhood of the origin by the power series
    \[
        \text{Li}_{n_1,\dots, n_r}(z_1,\dots, z_r)=\sum_{0<k_1<\dots <k_r} \frac{z_1^{k_1}\dots z_r^{k_r}}{k_1^{n_1}\dots k_r^{n_r}}  
    \]
    in the complex variables $z_1,\dots, z_r$. This power series converges absolutely for $|z_i|<1$. If $n_r\ge 2$, convergence extends to $|z_i|\le 1$. The quantity $n_1+\dots+n_r$ is called \emph{weight} of the multiple polylogarithm $\text{Li}_{n_1,\dots, n_r}(z_1,\dots,z_r)$.\\
    Given $w=[\omega_{i_s}\vert\dots \vert \omega_{i_0}]\in B$ with $\omega_{i_0}\neq 0$, we may write
    \[
        w=[\omega_0^{n_r-1}\vert \omega_{j_r}\vert \omega_0^{n_{r-1}-1}\vert \omega_{j_{r-1}}\vert\dots \vert \omega_0^{n_1-1}\vert \omega_{j_1}]
    \]
    with $\omega_{j_1}=\omega_{i_0}\neq 0$. For $k=1,\dots, n$, set $h_k=n-j_k+1$ and define $y_1,\dots, y_n$ as follows:
    \[
        y_1=(\sigma_2\dots \sigma_n)^{-1},\;\dots,y_{n-2}=(\sigma_2\sigma_3)^{-1}, \;y_{n-1}=\sigma_2^{-1}.
    \]
    Thus, we have
    \[
        \sigma_2=y_{n-1}^{-1},\; \sigma_3=(y_{n-2}y_{n-1})^{-1}, \; \dots,\; \sigma_l=(y_1\dots y_{n-1})^{-1}.
    \]
    By computing explicitly the derivatives of the polylogarithms using their series expansion, it is possible to see that 
    \[
        L_w(y_n)= (-1)^r\text{Li}_{n_1,\dots, n_r}\left( \frac{y_{h_1}\dots y_{n}}{y_{h_2}\dots y_n},\dots, \frac{y_{h_{r-1}}\dots y_n}{y_{j_r}\dots y_n},y_{h_r}\dots y_n\right).
    \]
    Consider now an algebraic differential form $\omega$ over $\f{M}_{0,l+3}$. By \cite[Section 8.3]{Brown-Multiple_zeta_values_and_periods_of_moduli_spaces_of_curves}, $\omega$ admits an analytic primitive with coefficients given by regular functions on $\f{M}_{0,l+3}$ and hyperlogarithms in the cubical coordinate $x_l$ in the algebra $B$ described above. The ponits $\sigma_1,\dots, \sigma_n$ have to be chosen as suitable products of the remaining cubical coordinates $x_1,\dots, x_{l-1}$ and their inverses. Instead of exposing a precise statement, we will turn to specific instances of this result later with the concrete example of the integrals $\xi_l$. \\
    As a final piece of notation, given a hyperlogarithm $L_w(z)$, we will use the symbol of indefinite integral
    \[
        \int L_w(z)dz
    \]
    to denote the unique primitive of $L_w(z)dz$ in $B$ whose regularization at $z=0$ takes the value $0$. Moreover, for multiple zeta values we follow the notational conventions of \cite{Burgos-Gil-Fresan-Multiple_zeta_values_from_numbers_to_motives}.\\

    \subsection{Examples in low dimensions}
        
    From now on we set $\Delta_l=[0,1]^l$. Let us give a few examples in low dimension of the general computation of $\xi_l$, starting with the case $l=2$. Although this specific integral is rather easy to compute, we illustrate the strategy step by step. For a fixed $x_2\in \C\setminus\{0\}$ we may regard the rational function $(1-x_1x_2)^{-1}$ as a hyperlogarithm in the variable $x_1$ over $\mathbb{P}^1(\C)\setminus \{0,1,x_2^{-1},\infty\}$. \\
    First, we exhibit a primitive of this function in the variable $x_1$ which is regularized in such a way that it vanishes for $x_1=0$, without logarithmic singularities. Such primitive is given by $-x_2^{-1}\log(1-x_1x_2)$. \\
    Next, we evaluate this primitive at $x_1=1$ to obtain $-x_2^{-1} \log(1-x_2)$. Thus,
    \[
        \int_{\Delta_2} \omega_2 = -\int_0^1\frac{\log(1-x_2)}{x_2} dx_2.
    \]
    With respect to the variable $x_2$, a regularized primitive for $-x_2^{-1}\log(1-x_2)$ is $L_{a_0a_1}(x_2)=\text{Li}_2(x_2)$, where
    \[
        a_0a_1=\left[ \frac{dt}{t} \,\middle\vert\, \frac{dt}{1-t} \right].
    \]
    We conclude that 
    \[
        \int_{\Delta_2} \omega_2 = \text{Li}_2(1) = \zeta(2).
    \]
    $ $\\
    
    Let us pass to $l=3$. In order to make our notation lighter, in the reduced bar complex of $\mathbb{P}^1(\C)\setminus\{0,\sigma_1,\dots, \sigma_l,\infty\}$ the element $[(\sigma_i-t)^{-1}dt]$ shall be denoted by $[\sigma_i]$, and similarly for all other words. We also write $[0]$ for $[t^{-1}dt]$. We need to deal with hyperlogarithms in one variable $x_j$ over $\mathbb{P}^1(\C)\setminus\{\sigma_0,\dots, \sigma_l,\infty\}$, where the $\sigma_i$'s depend on the variables $x_{j+1},\dots, x_l$. The first steps of the previous case imply that
    \[
        \int_{\Delta_3} \omega_3 = -\int_0^1 \int_0^1 \frac{\log(1-x_2)}{x_2(1-x_2x_3)}dx_2dx_3
    \]
    By decomposing into partial fractions, we see that
    \[
        \int_{\Delta_3} \omega_3 = -\int_0^1 \int_0^1 \frac{\log(1-x_2)}{x_2}dx_2dx_3+\int_0^1\int_0^1\frac{-x_3\log(1-x_2)}{1-x_2x_3}dx_2dx_3.
    \]
    The value of the first term is $\zeta(2)$. A primitive of the second integrand in the variable $x_2$ is $L_{[x_3^{-1}\vert  1]}(x_2)$. As a result,
    \[
        \int_{\Delta_3} \omega_3 = \zeta(2) +\int_0^1 L_{[x_3^{-1}\vert  1]}(x_2)\vert_{x_2=1} dx_3
    \]
    In order to find a primitive of $ L_{[x_3^{-1}\vert  1]}(x_2)\vert_{x_2=1}$ in the variable $x_3$, we first need to write this function as a hyperlogarithm on $\mathbb{P}^1(\C)\setminus\{0,1,\infty\}$. \\
    We achieve this as follows: first, we express $L_{[x_3^{-1}\vert  1]}(x_2)$ as a hyperlogarithm in $x_3$ over $\mathbb{P}^1(\C)\setminus\{0,1,x_2^{-1}\infty\}$; after this, we set $x_2=1$, which yields a hyperlogarithm over $\mathbb{P}^1(\C)\setminus\{0,1,\infty\}$. The first step requires the equality
    \[
        L_{[x_3^{-1}\vert  1]}(x_2)=L_{[0\vert x_2^{-1}]}(x_3)+L_{[1\vert x_2^{-1}]}(x_3)+(L_{[ 1]}(x_3)-L_{[x_2^{-1}]}(x_3))L_{[1]}(x_2).
    \]
    We will explain how to derive this equality later in Lemma~\ref{lemma: swap two variables in a hyperlog}. For the moment, this equation can be checked by taking the derivative with respect to $x_3$ and arguing that both sides vanish at the origin.\\
    Setting $x_2=1$, we infer that
    \[
        L_{[x_3^{-1}\vert  1]}(x_2)\vert_{x_2=1}= L_{[0\vert1]}(x_3)+L_{[1\vert 1]}(x_3).
    \]
    Hence,
    \[
        \int_{\Delta_3}\omega_3 = \zeta(2)+\int_0^1\left(L_{[0\vert 1]}(x_3)+L_{[1\vert 1]}(x_3) \right)dx_3
    \]
    Once again, we need to find a primitive for the integrand on the right-hand side in the variable $x_3$. Integration by parts yields 
    \begin{gather*}
        \int  L_{[0\vert 1]}(x_3) dx_3 = x_3 L_{[0\vert 1]}(x_3)-\int L_{[1]}(x_3)dx_3, \\
        \int  L_{[1\vert 1]}(x_3) dx_3 = (x_3-1) L_{[0\vert 1]}(x_3)+\int L_{[1]}(x_3)dx_3.
    \end{gather*}
    Overall, evaluating these primitives at $x_3=1$ finally gives
    \[
        \int_{\Delta_3}\omega_3 = \zeta(2)+\text{Li}_2(1)=2\zeta(2).
    \]
    \begin{remark}
        As is already apparent from this example, the primitive of $\omega_3$ produced by this algorithm is not algebraic, even if we have seen that $\omega_3$ is exact as an algebraic form. The use of hyperlogarithms, despite enabling precise computations of period integrals over $\overline{\f{M}}_{0,l+3}$, does not allow us to see which weights appear a priori.
    \end{remark}

    The case $l=4$ becomes rather complicated to handle in detail. Thus, we make explicit the pattern of partial fraction decompositions that has appeared so far. As a consequence, we reduce the computation of $\xi_l$ to simpler integrals.\\
    Write $f_l=(1-x_1x_2)\dots (1-x_{l-1}x_l)$, so 
    \[
    \int_{\Delta_l}\omega_l=\int_{\Delta_l} \frac{1}{f_l} dx_1\dots dx_l.
    \]
    Let $F^{(l)}_k(x_k,x_{k+1})$ be the $k$-th primitive found in the algorithm. This is defined inductively to be a primitive of $(1-x_kx_{k+1})^{-1}F^{(l)}_{k-1}(1,x_k)$, regularized to vanish at the origin. This means that for $k=1,\dots, l-2$
    \[
        \int_{\Delta_l}\omega_l =\int_{[0,1]^{l-k}} \frac{F^{(l)}_k(1,x_{k+1})}{(1-x_{k+1}x_k)\dots (1-x_{l-1}x_l)}dx_{k+1}\dots dx_l=\int_0^1 F^{(l)}_{l-1}(1,x_l)dx_l.
    \]
    Note that 
    \begin{align*}
        \int_{\Delta_l} \omega_l &= \int_0^1\int_0^1 \frac{1}{1-x_{l-1}x_l} \left( \int_{[0,1]^{l-2}} \frac{1}{f_{l-1}} dx_1\dots dx_{l-2}\right) dx_{l-1} dx_l \\
                                 &= \int_0^1\int_0^1 \frac{F^{(l-1)}_{l-2}(1,x_{l-1})}{1-x_{l-1}x_l} dx_{l-1} dx_l,
    \end{align*}
    which proves that $F^{(l)}_k=F^{(l-1)}_k$ for all $k=1,\dots, l-2$, by uniqueness of the primitives involved. Thus, we may focus on $F_l=F^{(l)}_{l-1}$, so 
    \[
        F_l(x_{l-1},x_l)=\int \frac{F_{l-1}(1,x_{l-1})}{1-x_{l-1}x_l}dx_{l-1}.
    \]
    Our previous computations give
    \begin{gather*}
        F_2(1,x_2)=-\frac{1}{x_2}\log(1-x_2); \\
        F_3(1,x_3)=\zeta(2)+L_{[0\vert 1]}(x_3)+L_{[1\vert 1]}(x_3).
    \end{gather*}
    The first step to compute $F_3$ was a partial fraction decomposition, which we now rewrite:
    \begin{align*}
        F_3(x_2,x_3)&=\int \frac{F_2(1,x_{2})}{1-x_{2}x_3}dx_{2} \\
        &= \int F_2(1,x_{2}) dx_2 + \int \frac{x_2x_3 F_2(1,x_{2})}{1-x_2x_3}dx_2.
    \end{align*}
    Let $G_3(x_2,x_3)$ be the second indefinite integral, which means that it is a primitive of $x_2x_3(1-x_2x_3)^{-1} F_2(1,x_{2})$ with respect to $x_2$. Evaluating at $x_2=1$ then yields
    \[
        F_3(1,x_3)= \int_{\Delta_2}\omega_2 + G_3(1,x_3).
    \]
    In the above expression, we have thus isolated $G_3(1,x_3)=L_{[0\vert 1]}(x_3)+L_{[1\vert 1]}(x_3)$. When integrating $F_3(1,x_3)$ on $[0,1]$ to compute $\int_{\Delta_3}\omega_3$, the first summand has weight $2$, and the same holds for the second one, because it is not necessary to increase the weight in order to find a primitive of $G_3(1,x_3)$. \\
    For $l=4$, no partial fraction decomposition is required. Let $G_4(x_3,x_4)$ be a primitive of $(1-x_3x_4)^{-1}G_3(1,x_3)$. Then
    \begin{align*}
        F_4(x_{3},x_4)&=\int \frac{F_3(1,x_3)}{1-x_{3}x_4}dx_{3}\\
        &= \int_{\Delta_2}\omega_2 \int \frac{1}{1-x_3x_4}dx_3 + \int \frac{G_3(1,x_3)}{1-x_3x_4}dx_3 \\
        &= \left(\int_{\Delta_2}\omega_2 \right) F^{(2)}_1(x_3,x_4) + G_4(x_3,x_4).
    \end{align*}

    Let us apply this strategy to the general case. Write $G_2=F_2$ and let $F_1$ be the constant function $1$. We define for $l\ge 2$
    \[
        G_l(x_{l-1},x_l)=
        \begin{dcases}
            \int \frac{G_{l-1}(1,x_{l-1})}{1-x_{l-1}x_l}dx_{l-1} & \text{if $l$ is even;}\\
            \int \frac{x_{l-1}x_lG_{l-1}(1,x_{l-1})}{1-x_{l-1}x_l}dx_{l-1} & \text{if $l$ is odd.}
        \end{dcases}
    \]
    \begin{lemma}
    \label{lemma: introduction of the coefficients a_lk}
    For $l\ge 2$, we have
        \begin{align*}
        F_l(1,x_l)=\sum_{\text{$2 \le k \le l-2$, $k$ even}} a_{l,k} F_k(1,x_l) + G_l(1,x_l) & \quad \text{if $l$ is even}; \\
        F_l(1,x_l)=\sum_{\text{$1 \le k \le l-2$, $k$ odd}} a_{l,k} F_k(1,x_l) + G_l(1,x_l) & \quad \text{if $l$ is odd}; 
        \end{align*}
        where the $a_{l,k}$'s are defined recursively as follows:
        \begin{align*}
            & a_{l,k}=a_{l-1,k-1} & \text{if $l\ge 3$, $k\ge 2$}; \\
            & a_{l,1}=\xi_{l-1}-\sum_{\text{$2\le k\le l-3$, $k$ even}}a_{l-1,k}\xi_k & \text{if $l$ is odd.} 
        \end{align*}
    \end{lemma}
    \begin{proof}
        We argue by induction, the cases $l=2,3,4$ being already verified. Suppose that $l$ is odd, so we prove the formula for $l+1$ even. We have
        \begin{align*}
            F_{l+1}(x_l,x_{l+1}) & = \int \frac{F_l(1,x_l)}{1-x_lx_{l+1}}dx_l\\
            & = \sum_{\text{$k$ odd}} a_{l,k} \int\frac{F_k(1,x_l)}{1-x_lx_{l+1}}dx_l + \int \frac{G_l(1,x_l)}{1-x_lx_{l+1}}dx_l \\
            &=\sum_{\text{$k$ odd}} a_{l,k} F_{k+1}(x_l,x_{l+1}) + G_{l+1}(x_l,x_{l+1})\\
            &= \sum_{\text{$2\le k \le l$, $k$ even}} a_{l,k-1} F_{k}(x_l,x_{l+1})+ G_{l+1}(x_l,x_{l+1}),
        \end{align*}
        which verifies the claim for $l+1$ even. Note that $a_{l+1,k}=a_{l,k-1}$.\\
        Suppose now that $l$ is even, in which case we decompose into partial fractions as follows:
        \begin{align*}
            F_{l+1}(x_l,x_{l+1}) & = \int \frac{F_l(1,x_l)}{1-x_lx_{l+1}}dx_l\\
            & = \sum_{\text{$k$ even}} a_{l,k} \int\frac{F_k(1,x_l)}{1-x_lx_{l+1}}dx_l + \int \frac{G_l(1,x_l)}{1-x_lx_{l+1}}dx_l \\
            & = \sum_{\text{$k$ even}} a_{l,k} F_{k+1}(x_l,x_{l+1})+\int G_l(1,x_l) dx_l + \int\frac{x_l x_{l+1} G_l(1,x_l)}{1-x_lx_{l+1}}dx_l\\
            & =\sum_{\text{$k$ even}} a_{l,k} F_{k+1}(x_l,x_{l+1})+\int G_l(1,x_l) dx_l + G_{l+1}(x_l,x_{l+1}).
        \end{align*}
        When setting $x_l=1$ the integral in the middle coincides with
        \begin{align*}
            \int_0^1 G_l(1,x_l)dx_l &= \int_0^1 F_l(1,x_l) dx_l - \sum_{\text{$2 \le k \le l-2$, $k$ even}} a_{l,k} \int_0^1 F_k(1,x_l) dx_l\\
            &= \int_{\Delta_l}\omega_l - \sum_{\text{$2 \le k \le l-2$, $k$ even}} a_{l,k}\int_{\Delta_k}\omega_k.
        \end{align*} 
        Defining this quantity as $a_{l+1,1}$ yields the statement.
    \end{proof}

    \begin{corollary}
    \label{cor: shape of a_lk}
        Giving degree $n$ to $\xi_n$, the coefficient $a_{l,k}$ is a homogeneous polynomial in the $\xi_n$'s of degree $l-k$. Moreover, only $\xi_n$'s with $n$ even appear.
    \end{corollary}
    \begin{proof}
        Since $a_{l,k}=a_{l-k+1,1}$, it is enough to prove the statement for $a_{l,1}$ with $l$ odd. However, we have 
        \[
            a_{l,1}=\xi_{l-1}-\sum_{\text{$2\le k\le l-3$, $k$ even}}a_{l-1,k}\xi_k,
        \]
        so the claim follows by induction.
    \end{proof}
    \noindent
    The main key step that remains to be addressed is to find an explicit expression for $G_l(1,x_l)$ for all $l\ge 2$. This will be the main focus of the remainder of the article.

    \subsection{A recurrence sequence}
        
    Before proceeding, we introduce a recurrence sequence of complex numbers which plays a crucial role in the computation of $G_l(1,x_l)$. Given the purely combinatorial arguments used to understand this sequence, we prefer to study it in detail in this separate section.\\
    Throughout this section, we fix a sequence of complex numbers $\psi_0,\psi_2,\psi_4,\dots$, which will later be set equal to the numbers $\psi_n$ defined in the introduction together with $\psi_0=1$. However, for the moment we allow ourselves a slightly more general setup by letting these $\psi_n$ assume arbitrary complex values.\\
    For all integers $l\ge 1$ and $1\le m\le l-1$ with $m\equiv l \;(\text{mod}\; 2)$, we define recursively the numbers
    \[
        \beta^{(l+1)}_m= \sum_{\substack{n=0 \\ \text{ $n$ even}}}^{l-m}  \psi_n\beta^{(l)}_{m-1+n},
    \]
    with initial conditions $\beta^{(l)}_0=0$ and $\beta^{(l+1)}_l=1$ for all $l\ge 2$.\\
    We start by giving a more explicit description of these numbers.
    \begin{lemma}
        \label{lemma: closed formula for the beta's}
        For $s=1,\dots, l-1$ let $K_{l+1,m,s}$ be the set of $s$-tuples 
        \[
            (k_1,\dots, k_s)\in \left\{0,\dots, \frac{l-m}{2}\right\}^{s}
        \]
        satisfying the following conditions:  
        \begin{enumerate}
            \item $2(k_1+\dots+k_s)=l-m$;
            \item $k_s\neq 0$;
            \item for all $r=1,\dots, s$ we have $m+\sum_{i=1}^r \left( 2k_i-1\right)\ge 1$.
        \end{enumerate}
        Let us also set 
        \[
            K_{l+1,m}=\bigcup_{s=1}^{l-1} K_{l+1,m,s}.
        \]
        Then the following formula holds:
        \[
            \beta^{(l+1)}_m=\sum_{(k_1,\dots,k_s)\in K_{l+1,m}} \psi_{2k_1}\dots \psi_{2k_s}.
        \]
    \end{lemma}
    \begin{proof}
        We argue by induction on $l$, so
        \[
            \beta^{(l+1)}_m= \sum_{n=0}^{\frac{l-m}{2}}  \psi_{2n}\beta^{(l)}_{m-1+n}= \psi_{l-m}+\sum_{n=0}^{\frac{l-m}{2}-1}  \psi_{2n} \sum_{(k_{n,1},\dots,k_{n,s})\in K_{l,m-1+2n}} \psi_{2k_{n,1}}\dots \psi_{2k_{n,s}}.
        \]
        In the last sum, if $m=1$ we agree that $K_{l,0}=\emptyset$ and the corresponding sum is zero. Given $n\in \{0,\dots, \frac{l-m}{2}-1\}$ and $(k_{n,1},\dots,k_{n,s})\in K_{l,m-1+2n}$, we wish to show that $(n,k_{n,1},\dots,k_{n,s})\in K_{l+1,m} $.\\
        Of course, $k_{n,s}\neq 0$. Moreover, 
        \[
            2(n+k_{n,1}+\dots + k_{n,s})=2n+l-1-m+1-2n=l-m.
        \]
        For all $r\ge 1$ we have
        \[
            m+2n-1+\sum_{i=1}^r (2k_i-1) \ge m+2n-1+2-m-2n\ge 1.
        \]
        We also have $m+2n-1\ge 1$, because the case $m=1$ excludes the summand corresponding to $n=0$. This shows that $(n,k_{n,1},\dots,k_{n,s})\in K_{l+1,m} $.\\
        Conversely, let $(h,h_1,\dots, h_s)\in K_{l+1,m}$. We then have $2(h_1+\dots+h_s)=l-m-2h$ and $h_s\neq 0$. For all $r=1,\dots, s$ we have $m+2h-1+2(h_1+\dots+h_r)-r\ge 1$. Thus, $(h_1,\dots, h_s)\in K_{l,m-1+2h}$. We conclude that all tuples in $K_{l+1,m}$ appear exactly once in the sum above, whence the statement.    
    \end{proof}
    \noindent
    To give some intuition to the formula in the previous lemma, consider a matrix which has $\beta^{(l)}_m$ in the $(l,m)$-entry. We construct a directed graph with vertices the entries of this matrix and edges given as follows. Fix $l,m$ with $l\equiv m\;(\text{mod}\; 2)$. If $m\ge 2$ the edges with starting point the $(l+1,m)$-entry are the set $\{ E_0,E_2,\dots, E_{l-m}\}$, where $E_k$ for $k=0,\dots, l-m$ even has endpoint the $(l,k-1)$-entry. If $m=1$, we define the edges analogously, but omitting $E_0$.\\
    To compute $\beta^{(l+1)}_m$, one sums over all possible paths which start from the $(l+1,m)$-entry and reach any entry of the form $(l'+1,l')$ for some $l'\ge 1$. The summand corresponding to the path which is the composition of $E_{2k_1},\dots, E_{2k_s}$ equals the product $\psi_{2k_1}\dots \psi_{2k_s}$. \\
    It is easy to check that for any such choice of paths we have $2(k_1+\dots+k_s)=l-m$. Any number $k_i\in \{0,\dots, \frac{l-m}{2}\}$ may appear, with the only condition that the composition of paths considered does not go out of the matrix. After choosing $E_{2k_{r-1}}$, one may choose $k_0$ to be any number from $2$ to $\frac{l-m}{2}-k_1\dots-k_{r-1}$. The choice $r=0$ is allowed only if the endpoint of $E_{2k_{r-1}}$ is not in the first column. This is equivalent to requiring that 
    \[
        m+2\sum_{i=1}^{r}(k_i-1)\ge 1,
    \]
    which leads to the condition in the statement of the previous lemma. An example for the entry $(7,5)$ is visualized in Figure~\ref{figure}.\\
    
    \begin{figure}
    \begin{framed}
    \caption{Example of the computation of $\beta^{(7)}_5$.}
    \label{figure}
    \[
    \begin{tikzcd}
         & \beta^{(3)}_2 & & & &\\
        \beta^{(4)}_1 \ar[ru, "\psi_2"' {yshift=6pt, xshift=2pt}] & & \beta^{(4)}_3 &&&\\
        & \beta^{(5)}_2 \ar[lu, "\psi_0" {yshift=6pt, xshift=-2pt}] \ar[ru, "\psi_2"' {yshift=6pt, xshift=2pt}] & & \beta^{(5)}_4 && \\
        \beta^{(6)}_1 \ar[ru, "\psi_2"'{yshift=6pt, xshift=2pt}] \ar[rrru, "\psi_4" {xshift=40pt, yshift=11}, bend right=10] & & \beta^{(6)}_3 \ar[lu, "\psi_0" {yshift=6pt, xshift=-2pt}] \ar[ru, "\psi_2"' {yshift=6pt, xshift=2pt}] & & \beta^{(6)}_5 & \\
        & \beta^{(7)}_2 \ar[lu, "\psi_0" {yshift=6pt, xshift=-2pt}] \ar[ru, "\psi_2"' {yshift=6pt, xshift=2pt}] \ar[rrru, "\psi_4" {xshift=40pt, yshift=11}, bend right=10] & & \beta^{(7)}_4 & & \beta^{(7)}_6
    \end{tikzcd}
    \]
    \end{framed}
    \end{figure}
    
    From now on we will assume $\psi_0=1$. In this case, the product $\psi_{2k_1}\dots \psi_{2k_s}$ does not depend on the number of edges $E_0$ among $E_{2k_1},\dots, E_{2k_s}$. Thus, we may naturally express $\beta^{(l+1)}_m$ as a sum over the partitions of $\frac{l-m}{2}$ with certain coefficients that keep track of the maximal number of occurrences of $E_0$. To make this explicit, we introduce a new piece of notation.\\
    Fix non-negative integers $a_1\le a_2 \le\dots \le  a_s$ and let $N(a_1,\dots,a_s)$ be the number of $s$-tuples $(r_1,\dots, r_s)\in \Z_{\ge 0}^s$ such that
    \[
        r_1\le a_1, \quad r_1+r_2\le a_2, \quad \dots, \quad r_1+\dots + r_s \le a_s. 
    \]
    Notice that, setting
    \[
        y_1=r_1, \quad y_2=r_1+r_2, \quad \dots, \quad y_s=r_1+\dots + r_s,
    \]
    the number $N(a_1,\dots, a_s)$ coincides with the number of non-decreasing sequences of the form $0\le y_1\le \dots \le y_s $ such that $y_i\le a_i$ for all $i=1,\dots, s$. \\
    First, suppose that $a_1,\dots,a_s$ all equal a fixed integer $a\ge 0$ and write for brevity $N_s(a)=N(a_1,\dots, a_s)$. We then have
    \[
        N_s(a)=\binom{a+s}{s},
    \]
    which follows in an elementary way from a stars and bars argument.\\
    We turn now to $N(a_1,\dots, a_s)$ for general $a_1,\dots, a_s$. For convenience of notation, for all $a\in \Z$ we set $N_0(a)=1$, while for $s\ge 1$ and $a<0$ we define $N_s(a)=0$. Also, the symbol $N_m(a_1,\dots, a_{n-1};a_n)$ will be used for $N(a_1,\dots, a_s)$ with $a_n=a_{n+1}=\dots=a_s$.
    \begin{lemma}
        \label{lemma: explicit formula for N(a1,...as)}
        The following formula holds:
        \[
            N(a_1,\dots, a_s)=\sum_{q_0=1}^s \sum_{\substack{q_1,\dots, q_{s-1}=0 \\ q_1+\dots +q_{s-1}=s-q_0 \\ q_{s-1}+\dots + q_{s-j}\le j }}^{s-q_0} \binom{a_1+q_0}{q_0} \prod_{j=1}^{s-1} \binom{a_{j+1}-a_j-1+q_j}{q_j}.
        \]
    \end{lemma}
    \begin{proof}
        We argue by induction on $s$. For $s=1$ and any $a_1\ge 1$ the induction basis is verified because $N(a_1)=\binom{a_1+1}{1}=a_1+1$. We assume that the claim holds for $s-1$ and any choice of $a_1,\dots, a_{s-1}$. \\
        The set of $y_1,\dots, y_s\in \Z$ such that $0\le y_1\le \dots \le y_s$ and $y_i\le a_i$ can be written as the disjoint union of the following sets. First, we consider the $s$-tuples $(y_1,\dots, y_s)$ among these for which $y_s\le a_1$, which are $N_s(a_1)$ in total. \\
        Next, we consider those that satisfy $y_{s-1}\le a_1$ and $a_1+1\le y_s \le a_2$: the choices for $(y_1,\dots, y_{s-1})$ are $N_{s-1}(a_1)$, while the ones for $y_s$ are $N(a_2-a_1-1)=a_2-a_1$. Proceeding like this, for all $m=1,\dots, s-1$ we isolate the $s$-tuples $(y_1,\dots, y_s)$ which satisfy $y_m\le a_1$ and $a_1+1\le y_{m+1}\le \dots \le y_s \le a_2$, which amount to $N_m(a_1)N_{s-m}(a_2-a_1-1)$.\\
        After these,  for all $m=2,\dots, s-1$ we consider all $s$-tuples $(y_1,\dots, y_s)$ for which we have $y_{m}\le a_2$ and $a_2+1\le y_{m+1}\le \dots \le y_m \le a_3$. These are $N_{m-1}(a_1;a_2)N_{s-m}(a_3-a_2-1)$. Continuing in this manner, one is led to the following equality:
        \begin{equation*}
            N(a_1,\dots, a_s)=N_s(a_1)+\sum_{m=1}^{s-1} \sum_{n=1}^m N_{m-n+1}(a_1,\dots, a_{n-1};a_n) N_{s-m}(a_{n+1}-a_n-1).
        \end{equation*}
        For $m$ and $n$ in the above ranges, set $(a_1',\dots, a_n')=(a_1,\dots, a_{n-1},a_n,\dots, a_n)$. By induction hypothesis, we have
        \[
            N_{m-n+1}(a_1,\dots, a_{n-1};a_n)=\sum_{k_0=1}^m \sum_{k_1,\dots, k_{m-1}} N_{k_0}(a_1') \prod_{j=1}^{m-1}N_{k_j}(a'_{j+1}-a'_j-1),
        \]
        the internal sum running over all $k_1,\dots, k_{m-1}\in \{0,\dots, m-k_0\}$ that satisfy $k_1+\dots + k_{m-1}=m-k_0$ and $k_{m-1}+\dots +k_{m-j}\le j$. If $k_j\ge 1$ for $j\ge n$, then $N_{k_j}(a'_{j+1}-a'_j-1)=N_{k_j}(a_n-a_n-1)=N_{k_j}(-1)=0$. As a result, we may assume that $k_{n},\dots,k_{m-1}=0 $. Moreover, for $j\le n$ we have $a_j'=a_j$, so 
        \[
            N_{m-n+1}(a_1,\dots, a_{n-1};a_n)=\sum_{k_0=1}^m \sum_{k_1,\dots, k_{n-1}} N_{k_0}(a_1) \prod_{j=1}^{n-1}N_{k_j}(a_{j+1}-a_j-1).
        \]
        The number $N(a_1,\dots, a_s)$ is therefore the sum over $m=1,\dots, s$, $n=1,\dots, m$, $k_0=1,\dots, m$ and $k_1,\dots, k_{n-1}=0,\dots, m-k_0$ with $k_1+\dots + k_{n-1}=m-k_0$ and $k_{n-1}+\dots + k_n-j\le m-n+j$ of terms of the form
        \[
            N_{k_0}(a_1)\prod_{j=1}^{n-1}N_{k_j}(a_{j+1}-a_j-1) \, N_{s-m}(a_{n+1}-a_n-1).
        \]
        To rearrange this sum, let us define
        \[
            q_0=k_0, \quad q_i=k_i \; (i=1,\dots, n-1), \quad q_n=s-m, \quad q_j=0 \; (j=n+1,\dots, s-1).
        \]
        It follows immediately that $q_0=1,\dots, s$, every $q_i$ ranges from $0$ to $s-q_0$ and $q_1+\dots + q_{s-1}=s-q_0$. Moreover,
        \[
            q_{s-1}+\dots +q_{s-j}=q_n+(q_{n-1}+\dots + q_{n-(n-s+j)})\le s-q_0+m-n+(n-s+j)\le j.
        \]
        Given $(q_0,\dots, q_{s-1})$ as above, one can recover the quantities $m,n,k_0,\dots, k_{n-1}$ in an obvious way, choosing $n$ as the largest $j$ for which $q_j\ge 1$. The only case where this choice is not possible is when $q_0=s$, which yields the $(s-1)$-tuple $(0,\dots, 0)$. This bijection between the $q_j$'s and the $k_j$'s yields the equality
        \[
            N(a_1,\dots, a_s)=\sum_{q_0=1}^s \sum_{\substack{q_1,\dots, q_{s-1}=0 \\ q_1+\dots +q_{s-1}=s-q_0 \\ q_{s-1}+\dots + q_{s-j}\le j }}^{s-q_0} N_{q_0}(a_1) \prod_{j=1}^{s-1} N_{q_j}(a_{j+1}-a_j-1),
        \]
        which proves the statement.
    \end{proof}
    \begin{remark}
        The statement of this lemma is also valid when $a_1=0$. To have a more symmetric formula, we may set $a_0=0$ and observe that $N(a_0,a_1,\dots, a_s)=N(a_1,\dots, a_s)$, so
        \[
            N(a_1,\dots, a_s)=
            \sum_{\substack{q_1,\dots, q_{s}=0 \\ q_{s}+\dots + q_{s-j}\le j+1 }}^{s}
            \prod_{j=1}^{s} 
            \binom{a_{j}-a_{j-1}-1+q_j}{q_j}.
        \]
        However, we have preferred to write the formula for $N(a_1,\dots, a_s)$ as in the statement of the previous lemma because it leads to a more natural expression for the numbers $\gamma_{k_1,\dots,k_s}$ appearing in the next corollary.
    \end{remark}
    \begin{remark}
        The formula in Lemma~\ref{lemma: explicit formula for N(a1,...as)} is not meaningful for the case $a_i=i-1$ for $i=1,\dots, s$. To deal with this, we observe that 
        \[
            N(0,1,\dots, s-1)=\frac{1}{s+1}\binom{2s}{s}
        \]
        is the $s$-th Catalan number, which counts Dyck paths on a $s\times s$ grid \cite{Chung-Feller-On_fluctuations_in_coin-tossing} or the number of triangulations of a $(s+2)$-gon \cite{Crepinsek-Mernik-An_efficient_representation_for_solving_Catalan_number_related_problems}.
    \end{remark}
    \begin{corollary}
        \label{cor: formula for beta's without psi_0}
        In the notation of Lemma~\ref{lemma: closed formula for the beta's}, suppose that $\psi_0=1$. Let $P_{l+1,m}$ be the set of partitions of $\frac{l-m}{2}$, that is, the set of $s$-tuples $(k_1,\dots,k_s)$ of positive integers whose sum equals $\frac{l-m}{2}$. Given $(k_1,\dots,k_s)\in P_{l+1,m}$, let 
        \[
            \gamma_{k_1,\dots,k_s}= \sum_{q_0=1}^s \sum_{\substack{q_1,\dots, q_{s-1}=0 \\ q_1+\dots +q_{s-1}=s-q_0 \\ q_{s-1}+\dots + q_{s-j}\le j }}^{s-q_0} \binom{m-1+q_0}{q_0} \prod_{j=1}^{s-1} \binom{2k_j-2+q_j}{q_j}.
        \]
        Then 
        \[
            \beta^{(l+1)}_m = \sum_{(k_1,\dots, k_s)\in P_{l+1,m}} \gamma_{k_1,\dots,k_s} \psi_{2k_1}\dots \psi_{2k_s}.
        \]
    \end{corollary}
    \begin{proof}
        The elements of $K_{l+1,m}$ are of the form
        \[
            (0^{r_1},k_1,\dots, 0^{r_s}, k_s)
        \]
        for some $(k_1,\dots, k_s)\in P_{l+1,m}$ and suitable $r_1,\dots, r_s\ge 1$. Since $\psi_0=1$, different values of the $r_i$'s yield the same summand in the formula in Lemma~\ref{lemma: closed formula for the beta's}, namely $\psi_{2k_1}\dots \psi_{2k_s}$. To compute the coefficient $\gamma_{k_1,\dots, k_s}$ we only need to find the maximal possible values of $r_1,\dots, r_s$. \\
        For every $q=1,\dots s$, the first $r_1+\dots+r_q+q-1$ components of $(0^{r_1},k_1,\dots, 0^{r_s}, k_s)$ must satisfy
        \[
            m+2 \sum_{i=1}^{q-1} k_i \ge r_1+\dots+r_q+q.
        \]
        These inequalities give the number of choices for $r_1,\dots, r_s$, which is precisely
        \[
        \gamma_{k_1,\dots,k_s}= N(m-1, m-2+2k_1, \dots, m-s+2(k_1+\dots+k_{s-1})).
        \]
        The statement finally follows from the explicit description of $N(a_1,\dots,a_s)$ given in Lemma~\ref{lemma: explicit formula for N(a1,...as)}.
    \end{proof}
    \begin{remark}
        Notice that the value of $\gamma_{k_1,\dots, k_s}$ depends on the ordered $s$-tuple $(k_1,\dots, k_s)$ and not just on the unordered partition $k_1,\dots, k_s$. For example, for $l-m=6$, the ordered partitions $3=1+2$ and $3=2+1$ lead to
        \[
            \gamma_{1,2}=\frac{m(m+3)}{2}, \qquad \gamma_{2,1}=\frac{m(m+7)}{2}.
        \]
    \end{remark}
    \noindent
    From the previous corollary, we may easily write down a few explicit expressions for $\beta^{(l+1)}_m$. First, we remark that, for $a_1\le a_2\le a_3$,
    \begin{align*}
        N(a_1) &=a_1+1; \\
        N(a_1,a_2)&= \frac{1}{2}(2a_2+2-a_1)(a_1+1); \\
        N(a_1,a_2,a_3)&= \frac{1}{6}(a_1+1)(a_1^2-4a_1-3a_1a_3-3a_2^2+3a_2+6a_2a_3+6a_3+6).
    \end{align*}
    Then, let $l\ge 1$ and $1\le m\le l-1$ with $m\equiv l \;(\text{mod}\; 2)$. Assume that $\psi_0=1$. For small values of $l-m$, we have:
    \begin{align*}
        \beta^{(l+1)}_{l-2}=&(l-2)\psi_2; \\
        \beta^{(l+1)}_{l-4}=&(l-4)\psi_4 +\frac{(l-1)(l-4)}{2}\psi_2^2\\
        \beta^{(l+1)}_{l-6} =& (l-6)\psi_6 + (l-6)(l-1)\psi_2\psi_4  + \frac{(l-6)(l-2)(l-1)}{6}\psi_2^3.
    \end{align*}
    We conclude this section by looking at a generating function for the $\beta^{(l)}_m$'s. We first need the following lemma:
    \begin{lemma}
    \label{lemma: technical property of the beta's}
         Let $l$, $m$ be positive integers, $l\ge 2$, $2\le m\le l$ and $m\equiv l+1 \;(\textnormal{mod}\; 2)$. Then for all $q=1,\dots, m-1$ we have
        \[ 
            \beta^{(l+2)}_{m}=\sum_{\substack{k=0\\ k\equiv q+1\;(\textnormal{mod}\; 2)}}^{l-m+q} \beta^{(k+2)}_q \beta^{(l+1-k)}_{m-q}.
        \]
    \end{lemma}
    \begin{proof}
        We argue by induction on $l$. From the recursive definition of $\beta^{(l+2)}_m$ we deduce that
        \[
            \beta^{(l+2)}_m=\sum_{\substack{n=0\\ \text{$n$ even}}}^{l+1-m} \psi_n \beta^{(l+1)}_{m-1+n}.
        \]
        We may apply the induction hypothesis to the terms $\beta^{(l+1)}_{m-1+n}$ to obtain
        \begin{align*}
            \beta^{(l+2)}_m &= \sum_{\substack{n=0\\ \text{$n$ even}}}^{l+1-m} \psi_n  \sum_{\substack{k=2\\ k\equiv q+1\;(\textnormal{mod}\; 2)}}^{l-m-n+q} \beta^{(k+2)}_q \beta^{(l-k)}_{m+n-1-q} \\
            &= \sum_{\substack{k=0\\ k\equiv q+1\;(\textnormal{mod}\; 2)}}^{l-m+q}\beta^{(k+2)}_q \left( \sum_{\substack{n=0\\ \text{$n$ even}}}^{l-m-k+q} \psi_n \beta^{(l-k)}_{m+n-1-q}  \right).
        \end{align*}
        The last sum in brackets coincides with $ \beta^{(l-k+1)}_{m-q}$, whence the statement follows.
    \end{proof}
    \noindent 
    For all $m\ge 1$ we introduce the formal power series
    \[
        B_m = \sum_{n=0}^\infty \beta^{(m+1+2n)}_m t^n.
    \]
    Lemma~\ref{lemma: technical property of the beta's} implies that $B_m=B_1^m$. Thus, in view of the recursive definition of $\beta^{(l)}_m$, it follows that $B_1$ satisfies the equation
    \[
        B_1=1+ \psi_2B_1^2\,t+ \psi_4B_1^4 \, t^2 + \dots = \sum_{n=0}^\infty \psi_{2n}B_1^{2n}\,t^n.
    \]
    If we write
    \[
        \Psi=\sum_{n=0}^\infty \psi_{2n}t^n,    
    \]
    this argument proves the following
    \begin{corollary}
        \label{cor: generating series}
        The formal power series $B_1$ satisfies the equation
        \[
            B_1=\Psi(B_1^2\,t).
        \]
    \end{corollary}

    \subsection{Final computation}
        
    To find an explicit expression for $G_{l}(1,x_l)$ and thus compute $\xi_l$, we need to introduce some notation. The reduced bar complex of $\mathbb{P}^1(\C)\setminus \{0,1,\infty\}$ is isomorphic to the free shuffle algebra over two letters. To simplify notation, we denote these letters by $0$ and $1$, which correspond to the differential forms $dz/z$ and $dz/(1-z)$ respectively. A standard monomial in this algebra is therefore represented by a tuple of zeroes and ones in the bar notation: $i=\left[ 0^{n_{r}} \,\middle\vert\, 1 \,\middle\vert\, \dots \,\middle\vert\, 0^{n_1-1} \,\middle\vert\, 1\right]$ for some $n_1,\dots, n_r\ge 1$. The $k$-th component of $i$ will be denoted by $i_k$. Moreover, we write $w(i)$ for the \emph{weight} of $i$, i.e. the number of components of $i$, and $l(i)$ for the \emph{length} of $i$, i.e. the number of non-zero components of $i$. Given $\epsilon\in \{ 0,1\}^{r-1}$, $\epsilon=(\epsilon_1,\dots, \epsilon_{r-1})$, we set $i(\epsilon)=\left[0^{n_{r-1}-1}\,\middle\vert\, \epsilon_{r-1} \,\middle\vert\, \dots \,\middle\vert\, 0^{n_1-1} \,\middle\vert\, \epsilon_1\right]$.\\
    One of the key points in the computation of the integrals $\xi_l$ is to find $G_{l+1}(1,x_{l+1})$ starting from $G_l(1,x_l)$. This involves finding a primitive of certain hyperlogarithms in $x_l$ with poles at $x_{l+1}$ and restricting them to $x_{l}=1$. This last step does not allow us to have an immediate representation of $G_{l+1}(1,x_{l+1})$ as a hyperlogarithm in $x_{l+1}$. We circumvent this obstacle by making use of the following lemma. The strategy of the proof follows Panzer's algorithm \cite[Section 2.4]{Panzer-Algortihms_for_the_symbolic_integration_of_hyperlogarithms_with_applications_to_Feynman_integrals}.
    \begin{lemma}
        \label{lemma: swap two variables in a hyperlog}
        Let $i$ be a standard monomial in the letters $\{0,1\}$ in the reduced bar complex of $\mathbb{P}^1(\C)\setminus \{0,1,\infty\}$. Then the following formula holds:
        \begin{align*}
             L_{\left[x^{-1} \,\middle\vert\, i \right]}(y) &= \sum_{\left[ a\,\middle\vert\, b \right]=[i]} (-1)^{w(a)-l(a)} \sum_{\epsilon\in \{0,1\}^{l(a)}} L_{\left[ a(\epsilon)\,\middle\vert\, y^{-1} \right]}(x)  L_{\left[b \right]}(y) \\
             & \qquad -\sum_{\left[ a\,\middle\vert\, b \right]=[i], \, b_1=1} (-1)^{w(a)-l(a)}\sum_{\epsilon\in \{0,1\}^{l(a)}} L_{\left[ a(\epsilon)\,\middle\vert\, 1 \right]}(x) L_{\left[b \right]}(y).
        \end{align*}
        In particular, we have
        \[
            L_{\left[x^{-1} \,\middle\vert\, i \right]}(y)\vert_{y=1}= \sum_{\left[ a\,\middle\vert\, b \right]=[i],\; b_1=0} (-1)^{w(a)-l(a)} \sum_{\epsilon\in \{0,1\}^{l(a)}} L_{\left[ a(\epsilon)\,\middle\vert\,1 \right]} (x) L_{\left[b \right]}(1).
        \]
    \end{lemma}
    \begin{proof}
        Write $i=\left[ 0^{n_{r-1}} \,\middle\vert\, 1 \,\middle\vert\, \dots \,\middle\vert\, 0^{n_1-1} \,\middle\vert\, 1\right]$. Let $l=n_1+\dots+n_{r-1}+1$ and set for short $m_1=n_1-1, \dots, m_{r-1}=n_{r-1}-1$. For all $k=1,\dots, l-1$ we let $i^k$ be the $(l-2)$-tuple obtained from $i$ by neglecting the $k$-th component and $i_{(k)}$ be the $(l-1-k)$-tuple obtained from $i$ by neglecting the first $k$ components. \\
        For all $k=1,\dots, l-1$ let $\sigma_{k+1}\in \{0,1\}$ be the $k$-th component of $i$ and $\sigma_1=x^{-1}$. Also set $\delta_k=0$ if $\sigma_k=0$, otherwise $\delta_k=1$. We then have
        \begin{align*}
            \partial_x L_{\left[x^{-1} \,\middle\vert\, i \right]}(y) =& \sum_{k=2}^{l-1} \frac{\partial_x(\sigma_k-\sigma_{k+1})}{\sigma_k-\sigma_{k+1}}
            \left((-1)^{\delta_{k+1}}L_{\left[x^{-1} \,\middle\vert\, i^{k}\right]}(y)
            -(-1)^{\delta_k}L_{\left[x^{-1} \,\middle\vert\, i^{k-1} \right]}(y)\right)\\
            & + \frac{\partial_x(\sigma_1-\sigma_{2})}{\sigma_1-\sigma_{2}}
            \left((-1)^{\delta_{2}}L_{\left[x^{-1} \,\middle\vert\, i^{1}\right]}(y)+L_{\left[i \right]}(y)\right)\\
            & +\frac{\partial_x \sigma_1}{y-\sigma_1}L_{\left[ i\right]}(y)+\frac{\partial_x\sigma_l}{\sigma_l}L_{\left[x^{-1} \,\middle\vert\, i^{l-1} \right]}(y).
        \end{align*}
        Given that $\partial_x\sigma_k=0$ when $k\neq 1$, we have
        \begin{align*}
            \partial_x L_{\left[x^{-1} \,\middle\vert\, i \right]}(y)&= 
            \left( \frac{1}{x(\sigma_2x-1)}-\frac{1}{x(xy-1)} \right)L_{\left[ i\right]}(y)+\frac{(-1)^{\delta_2}}{x(\sigma_2x-1)}L_{\left[x^{-1} \,\middle\vert\, i^{1}\right]}(y)\\
            &= \left( \frac{y}{1-xy}- \frac{\sigma_2}{1-\sigma_2x}\right)L_{\left[ i\right]}(y)+\frac{(-1)^{\delta_2}}{x(\sigma_2x-1)}L_{\left[x^{-1} \,\middle\vert\, i^{1}\right]}(y).
        \end{align*}
        If $\sigma_2=0$, then
        \[
            \partial_x \left( L_{\left[x^{-1} \,\middle\vert\, i \right]}(y)\right)= \partial_x\left(L_{[y^{-1}]}(x)L_{[i]}(y)\right)-\frac{1}{x}L_{\left[x^{-1} \,\middle\vert\, i^{1}\right]}(y).
        \]
        On the other hand, if $\sigma_2=1$, we have
        \[
             \partial_x \left( L_{\left[x^{-1} \,\middle\vert\, i \right]}(y)\right)=\partial_x\left( L_{[y^{-1}]-[1]}(x)L_{\left[ i\right]}(y)\right) + \frac{1}{x(1-x)}L_{\left[x^{-1} \,\middle\vert\, i^{1}\right]}(y).
        \]
        Applying the formula for $\sigma_2=0$ a total of $m_{r-1}=n_{r-1}-1$ times yields
        \begin{align*}
            \partial_x \left( L_{\left[x^{-1} \,\middle\vert\, i \right]}(y)\right)&=
            \partial_x \left( \sum_{k=0}^{m_{r-1}-1}(-1)^k L_{\left[0^k \,\middle\vert\, y^{-1}\right]}(x)L_{[i_{(k)}]}(y)\right)\\
            & \qquad + (-1)^{m_{r-1}}\frac{1}{x}\int\frac{1}{x} \int \dots \int \frac{1}{x}  L_{\left[x^{-1} \,\middle\vert\, i_{(m_{r-1})} \right]} \, dx.
        \end{align*}
        Since $i_{m_{r-1}}$ has a $1$ in the first position, we deduce that 
        \begin{align*}
            \partial_x \left( L_{\left[x^{-1} \,\middle\vert\, i \right]}(y)\right)&=
            \partial_x \left( \sum_{k=0}^{m_{r-1}}(-1)^k L_{\left[0^k \,\middle\vert\, y^{-1}\right]}(x)L_{[i_{(k)}]}(y)\right)\\
            & \qquad -(-1)^{m_{r-1}}\partial_x L_{\left[0^{m_{r-1}} \,\middle\vert\, 1\right]}(x)L_{[i_{(m_{r-1})}]}(y) \\
            & \qquad + (-1)^{m_{r-1}}\frac{1}{x}\int\dots \int \frac{1}{x} \int \left( \frac{1}{x}+\frac{1}{1-x} \right) L_{\left[x^{-1} \,\middle\vert\, i_{(m_{r-1}+1)} \right]} \, dx.
        \end{align*}
        Applying this argument inductively, we may deduce a formula for $L_{\left[x^{-1} \,\middle\vert\, i \right]}(y)$ as a hyperlogarithm in the variable $x$. Given $\epsilon\in \{ 0,1\}^{r-1}$, $\epsilon=(\epsilon_1,\dots, \epsilon_{r-1})$, we set 
        \[
            i(\epsilon)=\left[0^{n_{r-1}-1}\,\middle\vert\, \epsilon_{r-1} \,\middle\vert\, \dots \,\middle\vert\, 0^{n_1-1} \,\middle\vert\, \epsilon_1\right].
        \]
        Moreover, for all $k=0,\dots, l-1$ let $i^{(k)}$ be the $k$-tuple consisting of the first $k$ components of $i$. Define also $s_k$ to be the maximum $t\in \{0,\dots, l-1\}$ such that $n_t\le k$. We then have 
        \begin{align*}
             L_{\left[x^{-1} \,\middle\vert\, i \right]}(y) &= 
            \sum_{k=0}^{l-1} (-1)^{k-s_k} \sum_{\epsilon\in \{0,1\}^{s_k}} L_{\left[i^{(k)}(\epsilon) \,\middle\vert\, y^{-1} \right]} (x) L_{\left[i_{(k)}\right]}(y)\\
            & \qquad - \sum_{t=0}^{r-1} (-1)^{n_t-1-t} \sum_{\epsilon\in \{0,1\}^{n_t-1}} L_{\left[i^{(n_t-1)}(\epsilon) \,\middle\vert\, 1 \right]} (x) L_{\left[i_{(n_t-1)}\right]}(y).
        \end{align*}
        More precisely, our argument shows that the derivatives in $x$ of the right-hand and left-hand side coincide. The equality then follows from the fact that both sides extend holomorphically to the origin with value zero. The formulae in the statement follow at once.
    \end{proof}
    \noindent
    We may now gain a first insight into the shape of $G_{l}(1,x_l)$. To this extent, let $I_{m}$ be the set of $m$-tuples $i=(i_1,\dots, i_m)$ with $i_k\in\{0,1\}$ and $i_m=1$. Given $i\in I_m$, we also set 
    \[
            X(i)=\{ a\in \{0,1\}^{m-1} \mid \text{$a_k=1$ for all $k=1,\dots, m-1$ such that $i_k=1$ } \}.
    \]
    \begin{lemma}
        For every $l\ge 2$, we have
        \[
        G_l(1,x_l)=
        \begin{dcases}
            \frac{1}{x_l} \sum_{m=1}^{l-1}\sum_{i\in I_{m}} \alpha_{[i]}^{(l)} L_{[i]}(x_l) & \text{if $l$ is even;} \\
            \sum_{m=1}^{l-1}\sum_{i\in I_{m}} \alpha_{[i]}^{(l)} L_{[i]}(x_l) & \text{if $l$ is odd.}
        \end{dcases}
        \]
        Here, $\alpha_{[i]}^{(l)}$ is a linear combination with integer coefficients of multiple zeta values of weight exactly $l-1-w(i)$.
    \end{lemma}
    \begin{proof}
        We argue by induction on $l$. If $l$ is even, then 
        \begin{align*}
            G_{l+1}(x_l,x_{l+1})&=\int \frac{x_lx_{l+1}G_l(1,x_l)}{1-x_lx_{l+1}}dx_l = \sum_{m=1}^{l-1}\sum_{i\in I_{m}} \int \frac{\alpha_{[i]}^{(l)} L_{[i]}(x_l)}{x_{l+1}^{-1}-x_l}dx_l\\
            &=\sum_{m=1}^{l-1}\sum_{i\in I_{m}}\alpha_{[i]}^{(l)}L_{[x_{l+1}^{-1}\vert i]}(x_l).
        \end{align*}
        On the other hand, for $l$ odd we have
        \begin{align*}
            G_{l+1}(x_l,x_{l+1})&=\int \frac{G_l(1,x_l)}{1-x_lx_{l+1}}dx_l =\frac{1}{x_{l+1}} \sum_{m=1}^{l-1}\sum_{i\in I_{m}} \int \frac{\alpha_{[i]}^{(l)} L_{[i]}(x_l)}{x_{l+1}^{-1}-x_l}dx_l \\
            &=\frac{1}{x_{l+1}}\sum_{m=1}^{l-1}\sum_{i\in I_{m}} \alpha_{[i]}^{(l)} L_{[x_{l+1}^{-1}\vert i]}(x_l).
        \end{align*}
        Up to a possible factor $x_{l+1}^{-1}$, by Lemma~\ref{lemma: swap two variables in a hyperlog} $G_{l+1}(1,x_{l+1})$ coincides with
        \[
            \sum_{m=1}^{l-1}\sum_{i\in I_{m}} \alpha_{[i]}^{(l)} \sum_{\left[ a\,\middle\vert\, b \right]=[i],\; b_1=0} (-1)^{w(a)-l(a)} \sum_{\epsilon\in \{0,1\}^{l(a)}} L_{\left[ a(\epsilon)\,\middle\vert\,1 \right]} (x_{l+1}) L_{\left[b \right]}(1).
        \]
        We may then rewrite this expression for $G_{l+1}(1,x_{l+1})$ in the form presented in the statement by setting 
        \[
            \alpha_{[i]}^{(l+1)}= \sum_{n=2}^{l-m} \sum_{\substack{b\in I_{n} \\ b_1=0 }} L_{[b]}(1)
            \sum_{a\in X(i)} (-1)^{w(a)-l(a)} \alpha_{[a\,\vert\, b]}^{(l)} + \sum_{a\in X(i)\cap I_{m-1}} (-1)^{w(a)-l(a)} \alpha_{[a]}^{(l)}.
        \]
        Note that $I_0=\emptyset$. Also, observe that we are starting with the initial data $\alpha^{(2)}_\varnothing =0$, $\alpha_{[0]}^{(2)}=0$ and $\alpha_{[1]}^{(2)}=1$. The fact that the sum over $n$ starts at $n=2$ and not at $n=1$ depends on the fact that $I_1=\{[1]\}$ but $b_1=0$.
    \end{proof}
    \noindent
    For a more precise structure of the function $G_{l}(1,x_{l})$, we need the following definition:
    \begin{definition}
        Given an elementary word $i$ in the reduced bar complex, we say that $i$ is \emph{admissible} if $i_k=1$ for all $k\equiv w(i)\;(\text{mod}\; 2)$. Given $m\ge 1$, let $\mathscr{I}_m$ denote the set of admissible words of weight $m$.
    \end{definition}
    \noindent
    Continuing with the notation of the previous section, for all $n\ge 0$ let us define $\psi_0=1$ and $\psi_{2n}$ for $n\ge 1$ as the sum of all multiple zeta values of weight $2n$ with only even arguments, that is,
    \[
        \psi_{2n}=\sum_{ \substack{k_1+\dots +k_s =n \\ k_i\ge 1} } \zeta(2k_1,\dots, 2k_s).
    \]
    Sums of multiple zeta values of fixed weight, length and depth are generated by single zeta values by a formula due to Zagier and Ohno \cite{Ohno-Zagier-Multiple_zeta_values_of_fixed_weight_depth_and_height}. The case of similar sums with the restriction to even arguments has been treated by Hoffman \cite{Hoffman-On_multiple_zeta_values_of_even_arguments} and will be used at the very end of the proof.\\
    Given this sequence of complex numbers $\{\psi_{2n}\}_{n\ge 0}$, we have a corresponding set of numbers $\beta^{(l+1)}_m$ for all $l\ge 1$, $1\le m\le l$ with $m\equiv l \;(\text{mod}\; 2)$, defined as in the previous section.
    \begin{proposition}
        \label{prop: shape of G_l}
        The function $G_{l+1}(1,x_{l+1})$ takes the following form:
        \[
            G_{l+1}(1,x_{l+1})=
            \begin{dcases}
                \frac{1}{x_{l+1}} \sum_{m=1}^l \beta_m^{(l+1)} \sum_{i\in \mathscr{I}_m} L_{[i]}(x_{l+1}) & \text{if $l$ is odd;} \\
                \sum_{m=1}^l \beta_m^{(l+1)} \sum_{i\in \mathscr{I}_m} L_{[i]}(x_{l+1}) & \text{if $l$ is even.}
            \end{dcases}
        \]
        Here, we agree that $\beta^{(l+1)}_m=0$ if $m\not\equiv l \;(\textnormal{mod}\; 2)$.
    \end{proposition}
    \begin{proof}
        We check that $\alpha_{[i]}^{(l+1)}$ is zero when $i$ is not admissible and equals $\beta^{(l+1)}_{w(i)}$ otherwise. In particular, $\alpha_{[i]}^{(l+1)}$ depends only on $w(i)$ when $i$ is admissible.\\
        Suppose that $w(i)=l-n'$ for $n'\ge 1$. Since $\beta^{(l+1)}_{l-1}=1$, the computations done for $G_2(1,x_2)$ and $G_3(1,x_3)$ verify the claim for $l=1,2$. By induction, we assume that the claim holds for all $\alpha^{(l)}_{[j]}$ with $w(j)\ge w(i)$. We shall prove the statement for $\alpha^{(l+1)}_{[i]}$. We write $\alpha^{(l+1)}_{[i]}=P_{i,0}+P_{i,2}+\dots+ P_{i,n'}$, where
        \[
            P_{i,0}= \sum_{a\in X(i)\cap I_{l-n'-1}} (-1)^{w(a)-l(a)} \alpha_{[a]}^{(l)}
        \]
        and for $n=2,\dots, n'$
        \[
            P_{i,n} = \sum_{\substack{b\in I_{n} \\ b_1=0 }} L_{[b]}(1)
            \sum_{a\in X(i)} (-1)^{w(a)-l(a)} \alpha_{[a\,\vert\, b]}^{(l)}.
        \]
        If $n'=l-1$, so $w(i)=1$, then $P_{i,0}=0$. Otherwise, by induction hypothesis, since $l-w(a)=l-w(i)$, the sum runs only over the words $a$ which are admissible. Moreover, their value is independent of $a$, so 
        \[
            P_{i,0}= \beta_{l-n'-1}^{(l)}\sum_{a\in X(i)\cap \mathscr{I}_{l-n'-1}} (-1)^{w(a)-l(a)}.
        \]
        For $n\ge 2$ , the words $[a\vert b]$ appearing in the innermost sum of $P_{i,n}$ have weight $l-n'-1+n$. Since $n'+1-n<n'$, we may apply the induction hypothesis to $\alpha^{(l)}_{[a\vert b]}$. If $b$ is not admissible, then $[a\vert b]$ is also not admissible, so we get a zero contribution. We may therefore restrict the sum over $b\in I_n$ with $b_1=0$ to admissible $b$'s. However, $b$ has to start with $0$, hence $w(b)$ must be even for $b$ to be admissible. This shows that $P_{i,n}=0$ if $n$ is odd.\\
        Assume that $w(b)=n$ is even. We may also restrict the sum over $a\in X(i)$ to the case when $[a\vert b]$ is admissible. Since $w(b)$ is even, this is equivalent to asking for $a$ to be admissible. We get:
        \[
            P_{i,n} = \beta^{(l)}_{l-n'-1+n} \sum_{\substack{b\in \mathscr{I}_{n} \\ b_1=0 }} L_{[b]}(1) 
            \sum_{a\in X(i)\cap \mathscr{I}_{l-n'-1} } (-1)^{w(a)-l(a)}.
        \]
        Observe that for $n$ in the above range
        \[
            \sum_{\substack{b\in \mathscr{I}_{n} \\ b_1=0 }} L_{[b]}(1) =\psi_n.
        \]
        Indeed, the multiple zeta values that appear in the left-hand side are associated with all possible binary sequences of the form
        \[ 
            i=
            \left[\,
            0 \,\middle\vert\, 
            1 \,\middle\vert\, 
            \epsilon_1\,\middle\vert\,
            1 \,\middle\vert\,
            \epsilon_2 \,\middle\vert\,
            1 \,\middle\vert\,
            \dots \,\middle\vert\,
            \epsilon_{n-1} \,\middle\vert\,
            1 \, \right]
        \]
        with $\epsilon_i\in \{0,1\}$. In particular, these binary sequences are obtained by chaining blocks of the form $\left[\, 1 \,\middle\vert\, 1 \, \right]$ or $\left[\, 0 \,\middle\vert\, 1 \, \right] $, always starting with $\left[\, 0 \,\middle\vert\, 1 \, \right] $. Since these two blocks end with $1$, when we rewrite $i$ in terms of multi-indices, the block $\left[\, 1 \,\middle\vert\, 1 \, \right]$ turns into the multi-index $(1,1)$ and the block $\left[\, 0 \,\middle\vert\, 1 \, \right]$ into the multi-index $(2)$. Thus, $\psi_{n}$ is the sum of all multiple zeta values $\zeta(q_1,\dots, q_s)$ of weight $n$ such that $q_s=2$ and $(q_1,\dots, q_{s-1})$ is obtained by chaining arbitrary sequences of the form $(1,1)$ or $(2)$. \\
        By \cite[Corollary 1.139]{Burgos-Gil-Fresan-Multiple_zeta_values_from_numbers_to_motives}, a multiple zeta value $\zeta(n_1,\dots,n_r)$ such that $n_i\ge 2$ for all $i=1,\dots,r$ equals
            \[
                \zeta(\, \overset{n_r-2}{\overbrace{1,\dots, 1}},2,\dots, \overset{n_1-2}{\overbrace{1,\dots, 1}},2 \,).
            \]
        We conclude that the sum of $L_{[b]}(1)$ for all $b\in \mathscr{I}_{n}$ with $b_1=0$ equals the sum of all multiple zeta values of weight $n$ with only even arguments, which is precisely $\psi_n$. \\
        Overall, considering that $\psi_0=1$, it follows that 
        \begin{align*}
            \alpha^{(l+1)}_{[i]} &= \left( \beta_{l-n'-1}^{(l)} + \sum_{\substack{n=2 \\ \text{ $n$ even}}}^{n'}  \psi_n  \beta^{(l)}_{l-n'-1+n} \right)\sum_{a\in X(i)\cap \mathscr{I}_{l-n'-1} } (-1)^{w(a)-l(a)} \\
                            &= \beta^{(l+1)}_{l-n'} \sum_{a\in X(i)\cap \mathscr{I}_{l-n'-1} } (-1)^{w(a)-l(a)}.
        \end{align*}
        by induction hypothesis $\beta^{(l)}_{l-n'-1+n}=0$ if $n'-n$ is odd. Since only even $n$'s appear in the expression above for $\alpha^{(l)}_{[i]}$, this occurs when $n'$ is odd. As a result, $n'$ being odd implies that $\alpha^{(l+1)}_{[i]}=0$, just like $\beta_{l-n'}^{(l+1)}$ by definition.\\
        To prove the claim, we are left to see that
        \[
            \sum_{a\in X(i)\cap \mathscr{I}_{l-n'-1} } (-1)^{w(a)-l(a)}=
            \begin{cases}
                1 & \text{if $i$ is admissible}; \\
                0 & \text{if $i$ is not admissible}. 
            \end{cases}
        \]
        To check this, suppose first that $i$ is admissible. Given $a\in X(i)$, by the very definition of $X(i)$ we must have $a_k=1$ for all $k\equiv w(i)\;(\text{mod}\; 2)$. On the other hand, if $a\in \mathscr{I}_{w(i)-1}$, it must be that $a_k=1$ for all $k\equiv w(i)-1 \;(\text{mod}\; 2)$. Thus, each component of $a$ equals $1$, so $X(i)\cap \mathscr{I}_{w(i)-1}$ contains only one element and the above sum reduces to $1$. \\
        Suppose on the other hand that $i$ is not admissible. Let $q_1,\dots, q_s\in \{1,\dotsm, w(i)\}$ be  the components of $i$ such that for all $t=1,\dots, s$ we have $q_t\equiv w(i) \;(\text{mod} \; 2)$ and $i_{q_t}=0$. Notice that $s\ge 1$ by assumption. Any $a\in X(i)\cap \mathscr{I}_{w(i)-1}$ satisfies $a_k=1$ for $k\equiv w(i)-1\;(\text{mod} \; 2)$, so $a$ is determined by setting some components among $q_1,\dots, q_s$ to $1$ and leaving the other ones to $0$. It follows that, for all fixed $d\ge 0$, the number of $a\in X(i)\cap\mathscr{I}_{w(i)-1}$ such that $w(a)-l(a)=d$ equals $\binom{s}{d}$. Hence,
        \[
            \sum_{a\in X(i)\cap \mathscr{I}_{l-n'-1} } (-1)^{w(a)-l(a)}=\sum_{d=0}^s (-1)^{d}\binom{s}{d}=(1+(-1))^s=0,
        \]
        as desired.
    \end{proof}
    \noindent
    Recall that $\xi_{l+1}$ equals the sum of $\int_0^1 G_{l+1}(1,x_{l+1})dx_{l+1}$ together with a term that is computed recursively via the $\xi_k$'s with $k\le l$. Let us focus on the integral of $G(1,x_{l+1})$ on $[0,1]$.
    \begin{lemma}
        \label{lemma: integral of G_l}
        For all $l\ge 1$, we have
        \[
            \int_0^1 G_{l+1}(1,x_{l+1})\, dx_{l+1} =
            \begin{dcases}
                \sum_{\substack{m=1 \\ \text{$m$ odd}}}^l \beta^{(l+1)}_m \psi_{m+1} & \text{if $l$ is odd}; \\
                \beta^{(l+2)}_1 & \text{if $l$ is even}.
            \end{dcases}
        \]
    \end{lemma}
    \begin{proof}
        For $l\ge 1$ and $m\ge 1$, $m\equiv l \; (\text{mod}\; 2)$, let us set for short
        \[
            H_{l+1,m}(x_{l+1})= \sum_{i\in \mathscr{I}_m} L_{[i]}(x_{l+1}). 
        \]
        By Proposition~\ref{prop: shape of G_l}, recall that 
        \[
            G(1,x_{l+1})= 
            \begin{cases}
                \sum_{\substack{m=1 \\ \text{$m$ odd}}}^l \frac{1}{x_{l+1}} \beta^{(l+1)}_m H_{l+1,m}(x_{l+1}) & \text{if $l$ is odd}; \\
                \sum_{\substack{m=1 \\ \text{$m$ even}}}^l \beta^{(l+1)}_m H_{l+1,m}(x_{l+1}) & \text{if $l$ is even}.
            \end{cases}
        \]
        If $l$ is odd, we have 
        \[
            \int \frac{1}{x_{l+1}}H_{l+1,m}(1,x_{l+1})dx_{l+1}= \sum_{i\in\mathscr{I}_m} L_{[0\vert i]}(x_{l+1}).
        \]
        Since $m$ is odd, every $i\in \mathscr{I}_m$ starts with $1$. Evaluating this primitive at $x_{l+1}=1$ gives therefore $\psi_{m+1}$. This shows that for $l$ odd
        \[
            \int_0^1 G_{l+1}(1,x_{l+1}) \,dx_{l+1} = \sum_{\substack{m=1 \\ \text{$m$ odd}}}^l \beta^{(l+1)}_m \psi_{m+1}.
        \]
        Assume now that $l$ is even, so $m=2n$ is also even. For $i\in \mathscr{I}_m$, say $i=[\delta_1\vert 1\vert \dots \vert \delta_{n}\vert 1]$ for some $\delta_1,\dots, \delta_{n}\in \{0,1\}$, define $f_{\delta_1}(x_{l+1})$ as $x_{l+1}$ is $\delta_1=0$ and $1-x_{l+1}$ if $\delta_1=1$. Integration by parts gives
        \begin{align*}
            \int_0^1 L_{[i]}(x_{l+1})dx_{l+1} =& \left( (-1)^{\delta_1}f_{\delta_1}(x_{l+1})L_{[i]}(x_{l+1})\right) \Big\vert_0^1 \\
                & -(-1)^{\delta_1}\int_0^1 L_{[1\vert \delta_2\vert\dots\vert 1\vert \delta_{n}\vert 1]}(x_{l+1}) dx_{l+1}.
        \end{align*}
        The first term vanishes when $\delta_1=1$, while it equals $L_{[i]}(1)$ if $\delta_1=0$. As a result, we have
        \begin{gather*}
            \int_0^1 H_{l+1,m}(x_{l+1}) \,dx_{l+1} = \sum_{i\in \mathscr{I}_m} \int_0^1 L_{[i]}(x_{l+1}) \,dx_{l+1} \\
            = \sum_{\substack{i\in \mathscr{I}_m\\ i_1=0}}  L_{[i]}(1)    
            +\sum_{\delta_2,\dots,\delta_{n}}\left( \sum_{\delta_1\in\{0,1\}} (-1)^{\delta_1} \right) \int_0^1 L_{[1\vert \delta_2\vert\dots\vert \delta_{n}\vert 1]}(x_{l+1}) \,dx_{l+1} \\
            = \sum_{\substack{i\in \mathscr{I}_m\\ i_1=0}}  L_{[i]}(1) = \psi_{m}.
        \end{gather*}
        We conclude that for $l$ even
        \[
            \int_0^1 G_{l+1}(1,x_{l+1})\,dx_{l+1}= \sum_{\substack{m=1 \\ \text{$m$ even}}}^l \beta^{(l+1)}_m \psi_m.
        \]
        The right-hand side coincides with the recursive definition of $\beta^{(l+2)}_1$, which concludes the proof.
    \end{proof}
    \noindent
    Let us now turn to the final computation of the integrals in Theorem~\ref{theo: main}. As seen in Lemma~\ref{lemma: introduction of the coefficients a_lk}, we have 
    \[
        \xi_{l+1}=  \sum_{\substack{k=1 \\ k\equiv l+1 \;(\text{mod}\; 2)}}^{l-1} a_{l+1,k}\xi_k + \int_0^1G_{l+1}(1,x_{l+1})\,dx_{l+1}.
    \]
    In the proof of Lemma~\ref{lemma: introduction of the coefficients a_lk} it is clear that $a_{l+1,1}=\int_0^1G_{l}(1,x_l)\,dx_l$ for all $l$ even. Moreover, $a_{l+1,k}=a_{l-k+2,1}$ for every $l$ and every $k\equiv l+1 \;(\text{mod}\; 2)$. Thus,
    \begin{align*}
        \xi_{l+1} &=  \sum_{\substack{k=1 \\ k\equiv l+1 \;(\text{mod}\; 2)}}^{l-1} \xi_k \int_{0}^1 G_{l-k+1}(1,x_{l-k+1})\,dx_{l-k+1} + \int_0^1G_{l+1}(1,x_{l+1})\,dx_{l+1} \\
                  & = \sum_{\substack{k=1 \\ k\equiv l+1 \;(\text{mod}\; 2)}}^{l-1} \sum_{\substack{m=1 \\ \text{$m$ odd}}}^{l-k} \xi_k \beta^{(l+1-k)}_m\psi_{m+1} + \int_0^1G_{l+1}(1,x_{l+1})\,dx_{l+1}.
    \end{align*}
    \begin{proposition}
        \label{prop: xi as beta}
        For all $l\ge 2$, we have 
        \[
            \xi_l= 
            \begin{cases}
                \beta^{(l+2)}_1 & \text{if $l$ is even}; \\
                \beta^{(l+2)}_2 & \text{ if $l$ is odd}.
            \end{cases}
        \]
    \end{proposition}
    \begin{proof}
        Suppose that $l$ is odd. We prove the statement for $\xi_{l+1}$ by induction on $l$. From our previous arguments we have
        \[
            \xi_{l+1}=\sum_{\substack{k=1 \\ \text{$k$ even}}}^{l-1} \, \sum_{\substack{m=1 \\ \text{$m$ odd}}}^{l-k} \beta^{(k+2)}_1\beta^{(l+1-k)}_m\psi_{m+1} +\sum_{\substack{m=1 \\ \text{$m$ odd}}}^l \beta^{(l+1)}_m \psi_{m+1}.
        \]
        We wish to prove that this equals 
        \[
            \beta^{(l+3)}_1=  \sum_{\substack{m=0 \\ \text{$m$ even}}}^{l+1} \psi_m \beta^{(l+2)}_{m}.
        \]
        Rearranging the sum in $\xi_{l+1}$ by collecting the $\psi_m$'s, we see that
        \begin{align*}
            \xi_{l+1} &= \sum_{\substack{m=1 \\ \text{$m$ odd}}}^{l-2}\left(\sum_{\substack{k=1 \\ \text{$k$ even}}}^{l-m} \beta^{(k+2)}_1\beta^{(l+1-k)}_m+\beta^{(l+1)}_m \right) \psi_{m+1} + \beta^{(l+1)}_l\psi_{l+1} \\
             &= \sum_{\substack{m=2 \\ \text{$m$ even}}}^{l-1}\left(\sum_{\substack{k=0 \\ \text{$k$ even}}}^{l-m+1} \beta^{(k+2)}_1\beta^{(l+1-k)}_{m-1} \right) \psi_{m} + \beta^{(l+2)}_{l+1}\psi_{l+1},
        \end{align*}
        using the fact that $\beta^{(l+1)}_l=1=\beta^{(l+2)}_{l+1}$. The statement then follows if we prove that for all $l$ odd, $2\le m\le l-1$, $m$ even, we have
        \[ 
            \sum_{\substack{k=0 \\ \text{$k$ even}}}^{l-m+1} \beta^{(k+2)}_1\beta^{(l+1-k)}_{m-1}= \beta^{(l+2)}_{m},
        \]
        which is precisely the content of Lemma~\ref{lemma: technical property of the beta's} for $q=1$.\\
        Suppose now that $l$ is even. We have, in view of Lemma~\ref{lemma: integral of G_l},
        \[
            \xi_{l+1}=\sum_{\substack{k=1 \\ \text{$k$ odd}}}^{l-1} \,\sum_{\substack{m=1 \\ \text{$m$ odd}}}^{l-k} \xi_k \beta^{(l+1-k)}_m\psi_{m+1} + \beta^{(l+2)}_1.
        \]
        By induction hypothesis, it follows that 
        \begin{align*}
            \xi_{l+1} &= \sum_{\substack{k=1 \\ \text{$k$ odd}}}^{l-1} \,\sum_{\substack{m=1 \\ \text{$m$ odd}}}^{l-k} \beta^{(k+2)}_2\beta^{(l+1-k)}_m \psi_{m+1} + \beta^{(l+2)}_1 \\
                      &= \sum_{\substack{m=2 \\ \text{$m$ even}}}^{l} \left( \sum_{\substack{k=1 \\ \text{$k$ odd }}}^{l-m+1} \beta^{(k+2)}_2 \beta^{(l+1-k)}_{m-1} \right) \psi_m +\beta^{(l+2)}_1.
        \end{align*}
        By the recursive definition of $\beta^{(l+3)}_2$, the claim follows from the equality
        \[
            \sum_{\substack{k=1 \\ \text{$k$ odd }}}^{l-m+1} \beta^{(k+2)}_2 \beta^{(l+1-k)}_{m-1} = \beta^{(l+2)}_{m+1},
        \]
        which has been proved in Lemma~\ref{lemma: technical property of the beta's} with $q=2$.
    \end{proof}
    \noindent
    Let $l\ge 2$ with $l=2n$ if $l$ is even and $l=2n+1$ if $l$ is odd. From Corollary~\ref{cor: formula for beta's without psi_0}, we have 
    \[
        \xi_l=\sum_{\substack{1\le k_1,\dots, k_s\le n \\ k_1+\dots+k_s=n}} \gamma_{k_1,\dots, k_s} \psi_{k_1}\dots \psi_{k_s},
    \]
    where
    \[
            \gamma_{k_1,\dots,k_s}= 
            \begin{dcases}
                \sum_{\substack{q_1,\dots, q_{s-1}=0  \\ q_{s-1}+\dots + q_{s-j}\le j }}^{s-1}
                \prod_{j=1}^{s-1} \binom{2k_j-2+q_j}{q_j} & \text{if $l$ is even}; \\
                \sum_{q_0=1}^s \sum_{\substack{q_1,\dots, q_{s-1}=0 \\ q_1+\dots +q_{s-1}=s-q_0 \\ q_{s-1}+\dots + q_{s-j}\le j }}^{s-q_0} q_0 \prod_{j=1}^{s-1} \binom{2k_j-2+q_j}{q_j} & \text{if $l$ is odd}.
            \end{dcases}
    \]
    \begin{corollary}
        Let $l=2m$ be even. Then
        \[
            \xi_{l+1}=\sum_{h=0}^{m} \xi_{2h}\xi_{l-2h}.
        \]
    \end{corollary}
    \begin{proof}
        From the previous Proposition together with Lemma~\ref{lemma: technical property of the beta's} it follows that
        \[
            \xi_{l+1}= \beta^{(l+3)}_{2}=\sum_{\substack{k=2 \\ \text{$k$ even}}}^l \beta^{(k+2)}_1\beta^{(l+2-k)}_1= \sum_{\substack{k=2 \\ \text{$k$ even}}}^l \xi_k\xi_{l-k},
        \]
        hence the statement.
    \end{proof}
    \noindent
    To conclude the proof of Theorem~\ref{theo: main}, in analogy with the notation of the previous section, we have
    \[
        1+\sum_{n=1}^\infty \xi_{2n}\,t^n= \sum_{n=0}^\infty \beta_1^{(2+2n)}\,t^n=B_1(t)
    \]
    by Proposition~\ref{prop: xi as beta}. Consider now the following generating function for the $\psi_{2n}$'s:
    \[
        \Psi(t)=\sum_{n=0}^\infty \psi_{2n}\,t^n.
    \]
    For all $n\ge 1$ and $r\le n$ define $E(2n,r)$ to be the sum of all multiple zeta values of weight $2n$ and depth $r$ with only even arguments, that is,
    \[
        E(2n,r)=\sum_{ n_1+\dots+n_r = n } \zeta(2n_1,\dots,2n_r).
    \]
    By \cite{Hoffman-On_multiple_zeta_values_of_even_arguments}, these numbers admit the generating function
    \[
        F(t,s)= 1+ \sum_{1\le r\le n } E(2n,r)\, t^n  s^r= \frac{\sin \left( \pi \sqrt{t(1-s)}\right)}{\sin\left(\pi\sqrt{t}\right) \sqrt{1-s}}.
    \]
    Notice, in particular, that $E(2n,r)$ is a rational multiple of $\pi^{2n}$. Since $\psi_{2n}$ equals the sum of $E(2n,r)$ with $r$ that runs through $1,\dots, n$, we deduce that 
    \[
        \Psi(t)=\lim_{s\to 1} F(s,t)= \frac{\pi\sqrt{t}}{\sin\left(\pi\sqrt{t}\right)}.
    \]
    Owing to Corollary~\ref{cor: generating series}, $B_1(t)$ satisfies the equality
    \[
        B_1(t)=\Psi(B_1(t)^2\, t)= \frac{\pi B_1(t)\,\sqrt{t}}{\sin\left(\pi B_1(t)\,\sqrt{t}\right)},
    \]
    which leads to the identity
    \[
        B_1(t)=\frac{\arcsin\left(\pi\sqrt{t}\right)}{\pi\sqrt{t}},
    \]
    as claimed in Theorem~\ref{theo: main}.\\

    $ $\\
    We conclude with a comparison with other integrals that appear in the literature. We were pointed out that the integrals $\xi_l$ resemble the ones considered by Zlobin in \cite{Zlobin-Rhin_integrals}, which are of the form
    \[
        \int_{[0,1]^l} \prod_{j=1}^{l-1} \frac{x_j^{a_j-1}(1-x_j)^{b_j-a_j-1}}{(1-x_jx_l)^{c_j}}x_l^{a_l-1}(1-x_l)^{b_l-a_l-1} \, dx_1\dots dx_l
    \]
    for suitable integer parameters $a_j$, $b_j$ and $c_j$. To highlight the difference with the $\xi_l$'s, consider the subfamily
    \[
        I_l=\int_{[0,1]^l} \frac{x_l^{l-2}}{(1-x_1x_l)(1-x_2x_l)\dots (1-x_{l-1}x_l)} \, dx_1\dots dx_l.
    \]
    Let us sketch how Panzer's algorithm applies to the integrals $I_l$. We have
    \begin{align*}
        I_l&=\int_{[0,1]^l} \frac{x_l^{l-2}}{(1-x_1x_l)(1-x_2x_l)\dots (1-x_{l-1}x_l)} \, dx_1\dots dx_l\\
            &= \int_{[0,1]^l} \frac{x_l^{-1}}{(x_{l}^{-1}-x_1)\dots (x_l^{-1}-x_{l-1})} \,dx_1\dots dx_l.
    \end{align*}
    Since the variable $x_j$ for $j\neq l$ appears in only one factor in the denominator, applying Panzer's algorithm for the variables $x_1,\dots, x_{l-1}$ means to find at each step a primitive of $(x_{l-1}^{-1}-x_j)^{-1}$ and evaluate it at $x_j=1$. Thus,
    \begin{align*}
        I_l&=\int_{[0,1]^l} \frac{x_l^{-1}}{(x_{l}^{-1}-x_1)\dots (x_l^{-1}-x_{l-1})} \, dx_1\dots dx_l\\
            &= \int_0^1 \frac{L_{[x_l^{-1}]}(1)^{l-1}}{x_l}\, dx_l = \int_0^1 \frac{L_{[1]}(x_l)^{l-1}}{x_l}\, dx_l.
    \end{align*}
    For the last primitive, observe that $L_{[1]}(x_l)^{l-1}=L_{[1]^{\sh(l-1)}}(x_l)=(l-1)! L_{[1^{l-1}]}(x_l)$, so 
    \begin{align*}
        I_l &= \int_0^1 \frac{L_{[1]}(x_l)^{l-1}}{x_l}\, dx_l = (l-1)!\int_0^1 \frac{L_{[1^{l-1}]}(x_l)}{x_l}\, dx_l \\
        & =(l-1)!\,L_{[0\vert 1^{l-1}]}(1)=(l-1)!\,\zeta(l).
    \end{align*}
    The main computational advantage is that the first $l-1$ primitives in the algorithm can be found independently of each other. This property is also shared by the general form with arbitrary parameters, because the only products of variables that appear are of the form $x_jx_l$. In contrast to this situation, in the integrals $\xi_l$ the shape of each primitive strongly depends on the primitive found in the previous step. This makes the search for primitives increasingly more difficult.

    \bibliographystyle{alpha}
    \bibliography{Literatur.bib}
\end{document}